\documentclass[11pt,a4paper]{amsart} 
\calclayout
\makeatletter 
\def\serieslogo@{} 
\def\@setcopyright{} 
\makeatother

\usepackage{amssymb,enumerate}
\usepackage[all]{xy}
\SelectTips{eu}{} 

\theoremstyle{plain}

\newtheorem{theorem}{Theorem}[section]
\newtheorem*{Theorem}{Theorem}

\newtheorem{proposition}[theorem]{Proposition}

\newtheorem{lemma}[theorem]{Lemma}

\newtheorem{corollary}[theorem]{Corollary}
\newtheorem*{Corollary}{Corollary}

\newtheorem{itheorem}{Theorem}

\theoremstyle{definition}

\newenvironment{bfchunk}{\begin{chunk}\textbf}{\end{chunk}}
\newtheorem{chunk}[theorem]{}

\newtheorem{subchunk}{}[theorem]

\theoremstyle{remark}
\newtheorem*{Question}{Question}
\newtheorem{remark}[theorem]{Remark}

\newtheorem*{Remark}{Remark}

\numberwithin{equation}{theorem}

\newcommand{\ac}{\mathrm{ac}}
\newcommand{\tac}{\mathrm{tac}}
\newcommand{\qis}{{\simeq}}

\newcommand{\Add}{\operatorname{Add}\nolimits}

\newcommand{\auss}{\mathsf{T}}
\newcommand{\kauss}{\mathsf{G}}
\newcommand{\bass}{\mathsf{S}}
\newcommand{\kbass}{\mathsf{F}}

\newcommand{\Coker}{\operatorname{Coker}\nolimits}
\newcommand{\Coloc}{\operatorname{Coloc}\nolimits}
\newcommand{\Loc}{\operatorname{Loc}\nolimits}
\newcommand{\Thick}{\operatorname{Thick}\nolimits}
\newcommand{\holim}{\operatorname{holim\,}\nolimits}

\newcommand{\col}{\colon}
\newcommand{\Cone}{\operatorname{Cone}\nolimits}

\newcommand{\dstar}{D^\ast}

\newcommand{\Ext}{\operatorname{Ext}\nolimits}

\newcommand{\Flat}{\operatorname{Flat}\nolimits}
\newcommand{\Inj}{\operatorname{Inj}\nolimits}
\newcommand{\proj}{\operatorname{prj}\nolimits}
\newcommand{\Proj}{\operatorname{Prj}\nolimits}

\newcommand{\hh}[1]{H(#1)}
\newcommand{\HH}[2]{H^{#1}(#2)}

\newcommand{\Hom}{\operatorname{Hom}\nolimits}
\newcommand{\RHom}{\operatorname{{\mathbf R}Hom}\nolimits}

\newcommand{\can}{\mathsf{can}}
\newcommand{\inc}{\mathsf{inc}}
\newcommand{\id}{\operatorname{id}\nolimits}
\newcommand{\shift}{\operatorname{shift}\nolimits}
\renewcommand{\Im}{\operatorname{Im}\nolimits}
\newcommand{\op}{\mathrm{op}}
\newcommand{\Ker}{\operatorname{Ker}\nolimits}

\newcommand{\kproj}{\operatorname{{\bf K}_{prj}}}
\newcommand{\kinj}{\operatorname{{\bf K}_{inj}}}

\newcommand{\lto}{\longrightarrow}
\newcommand{\xlto}[1]{\stackrel{#1}\lto}
\newcommand{\xto}{\xrightarrow}

\newcommand{\comp}{\mathop{\raisebox{+.2ex}{\hbox{$\scriptstyle\circ$}}}}

\newcommand{\Mod}{\operatorname{Mod}\nolimits}
\newcommand{\pd}{\operatorname{pd}\nolimits}
\newcommand{\rank}{\operatorname{rank}}
\newcommand{\soc}{\operatorname{soc}\nolimits}

\newcommand{\un}{\un}
\newcommand{\wh}{\widehat}

\def\dd{\partial}
\def\Si{{\sf\Sigma}}

\def\A{{\mathcal A}}
\def\B{{\mathcal B}}
\def\C{{\mathcal C}}

\def\S{{\mathcal S}}
\def\T{{\mathcal T}}

\def\bbZ{\mathbb Z}

\def\sfa{\mathsf a}
\def\sfb{\mathsf b}

\def\sfi{\mathsf i}
\def\sfp{\mathsf p}
\def\sfq{\mathsf q}
\def\sfr{\mathsf r}
\def\sfs{\mathsf s}
\def\sft{\mathsf t}
\def\sfu{\mathsf u}
\def\sfv{\mathsf v}
\def\sfF{\mathsf F}
\def\sfG{\mathsf G}

\def\bfD{\mathbf D}
\def\bfK{\mathbf K}
\def\bfL{\mathbf L}

\newcommand{\eps}{\varepsilon}
\newcommand{\fm}{\mathfrak{m}}
\newcommand{\fp}{\mathfrak{p}}

\newcommand{\ges}{{\scriptscriptstyle\geqslant}}
\newcommand{\les}{{\scriptscriptstyle\leqslant}}

\newcommand{\quism}{quasi-isomorphism}
\newcommand{\quic}{quasi-isomorphic}

\begin{document}

\title[Acyclicity versus total acyclicity]{Acyclicity versus total acyclicity \\
  for complexes over noetherian rings}

\author{Srikanth Iyengar}
\address{Srikanth Iyengar \\ Department of Mathematics\\
  University of Nebraska\\ Lincoln, NE - 68588\\ U.S.A.}  \email{iyengar@math.unl.edu}

\author{Henning Krause}
\address{Henning Krause\\ Institut f\"ur Mathematik\\
  Universit\"at Paderborn\\ 33095 Paderborn\\ Germany.}
\email{hkrause@math.uni-paderborn.de}

%\date{\today}
%\dedicatory{}  

\thanks{S. I. was partly supported by NSF grant DMS 0442242}



\begin{abstract}
  It is proved that for a commutative noetherian ring with dualizing complex the homotopy
  category of projective modules is equivalent, as a triangulated category, to the
  homotopy category of injective modules. Restricted to compact objects, this statement is
  a reinterpretation of Grothendieck's duality theorem. Using this equivalence it is
  proved that the (Verdier) quotient of the category of acyclic complexes of projectives
  by its subcategory of totally acyclic complexes and the corresponding category
  consisting of injective modules are equivalent. A new characterization is provided for
  complexes in Auslander categories and in Bass categories of such rings.
\end{abstract}

\maketitle

\section*{Introduction}
Let $R$ be a commutative noetherian ring with a dualizing complex $D$; in this article,
this means, in particular, that $D$ is a bounded complex of injective $R$-modules; see
Section \ref{Dualizing complexes} for a detailed definition.  The starting point of the
work described below was a realization that $\bfK(\Proj R)$ and $\bfK(\Inj R)$, the
homotopy categories of complexes of projective $R$-modules and of injective $R$-modules,
respectively, are equivalent. This equivalence comes about as follows: $D$ consists of
injective modules and, $R$ being noetherian, direct sums of injectives are injective, so
$D\otimes_R- $ defines a functor from $\bfK(\Proj R)$ to $\bfK(\Inj R)$.  This functor
factors through $\bfK(\Flat R)$, the homotopy category of flat $R$-modules, and provides
the lower row in the following diagram:
\[
\xymatrix{ 
\bfK(\Proj R)\ar@<-1ex>[rr]_-{\inc}&& \bfK(\Flat R)\ar@<-1ex>[ll]_-{\sfq}
  \ar@<-1ex>[rr]_-{D\otimes_R -}&& \bfK(\Inj R)\ar@<-1ex>[ll]_-{\Hom_R(D,-)}}
\]
The triangulated structures on the homotopy categories are preserved by $\inc$ and
$D\otimes_R-$. The functors in the upper row of the diagram are the corresponding right
adjoints; the existence of $\sfq$ is proved in Proposition~\eqref{flat:projection}.
Theorem \eqref{foxby} then asserts:

\begin{itheorem}
\label{ifoxby}
The functor $D\otimes_R-\col \bfK(\Proj R)\to \bfK(\Inj R)$ is an equivalence of
triangulated categories, with quasi-inverse $\sfq\comp\Hom_R(D,-)$.
\end{itheorem}

This equivalence is closely related to, and may be viewed as an extension of,
Grothendieck's duality theorem for $\bfD^f(R)$, the derived category of complexes whose
homology is bounded and finitely generated. To see this connection, one has to consider
the classes of compact objects -- the definition is recalled in \eqref{compacts} -- in
$\bfK(\Proj R)$ and in $\bfK(\Inj R)$.  These classes fit into a commutative diagram of
functors:
\[
\xymatrix{
  \bfK^c(\Proj R)\ar@{->}[rr]^-{D\otimes_R -}&&\bfK^c(\Inj R)\\
  \bfD^f(R)\ar@{<-}[u]_{\wr}^{\mathsf P} \ar@{->}[rr]^-{\RHom_R(-,D)}
  &&\ar@{<-}[u]^{\wr}_{\mathsf I} \bfD^f(R)}
\]
The functor $\mathsf P$ is induced by the composite
\[
\bfK(\Proj R)\xto{\Hom_R(-,R)}\bfK(R)\xto{\can}\bfD(R)\,,
\]
and it is a theorem of J{\o}rgensen \cite{Jo} that $\mathsf P$ is an equivalence of
categories. The equivalence $\mathsf I$ is induced by the canonical functor
$\bfK(R)\to\bfD(R)$; see \cite{Kr}.  Given these descriptions it is not hard to verify
that $D\otimes_R-$ preserves compactness; this explains the top row of the diagram.  Now,
Theorem~I implies that $D\otimes_R-$ restricts to an equivalence between compact objects,
so the diagram above implies $\RHom_R(-,D)$ is an equivalence; this is one version of the
duality theorem; see Hartshorne \cite{Ha}.  Conversely, given that $\RHom_R(-,D)$ is an
equivalence, so is the top row of the diagram; this is the crux of the proof of Theorem
\ref{ifoxby}.

Theorem \ref{ifoxby} appears in Section \ref{Foxby equivalence}.  The relevant definitions
and the machinery used in the proof of this result, and in the rest of the paper, are
recalled in Sections \ref{Triangulated categories} and \ref{Homotopy categories}.  In the
remainder of the paper we develop Theorem~\eqref{foxby} in two directions. The first one
deals with the difference between the category of acyclic complexes in $\bfK(\Proj R)$,
denoted $\bfK_{\ac}(\Proj R)$, and its subcategory consisting of totally acyclic
complexes, denoted $\bfK_\tac(\Proj R)$. We consider also the injective counterparts.
Theorems~\eqref{gdefect:proj} and \eqref{gdefect:inj} are the main new results in this
context; here is an extract:

\begin{itheorem}
\label{igdefect}
The quotients $\bfK_\ac(\Proj R)/\bfK_\tac(\Proj R)$ and $\bfK_\ac(\Inj R)/\bfK_\tac(\Inj
R)$ are compactly generated, and there are, up to direct factors, equivalences
\begin{gather*}
\Thick(R,D)/\Thick(R)   \xlto{\sim}
   \big[\big(\bfK_\ac(\Proj R)/\bfK_\tac(\Proj R)\big)^c\big]^\op\\
\Thick(R,D)/\Thick(R)   \xlto{\sim}
    \big(\bfK_\ac(\Inj R)/\bfK_\tac(\Inj R)\big)^c\,.
\end{gather*}
\end{itheorem}

In this result, $\Thick(R,D)$ is the thick subcategory of $\bfD^f(R)$ generated by $R$ and
$D$, while $\Thick(R)$ is the thick subcategory generated by $R$; that is to say, the
subcategory of complexes of finite projective dimension. The quotient
$\Thick(R,D)/\Thick(R)$ is a subcategory of the category $\bfD^f(R)/\Thick(R)$, which is
sometimes referred to as the stable category of $R$. Since a dualizing complex has finite
projective dimension if and only if $R$ is Gorenstein, one corollary of the preceding
theorem is that $R$ is Gorenstein if and only if every acyclic complex of projectives is
totally acyclic, if and only if every acyclic complex of injectives is totally acyclic.

Theorem \ref{igdefect} draws attention to the category $\Thick(R,D)/\Thick(R)$ as a
measure of the failure of a ring $R$ from being Gorenstein.  Its role is thus analogous to
that of the full stable category with regards to regularity: $\bfD^f(R)/\Thick(R)$ is
trivial if and only if $R$ is regular. See \eqref{gdefect:cat} for another piece of
evidence that suggests that $\Thick(R,D)/\Thick(R)$ is an object worth investigating
further.

In Section \ref{An example} we illustrate the results from Section \ref{Acyclicity} on
local rings whose maximal ideal is square-zero. Their properties are of interest also from
the point of view of Tate cohomology; see \eqref{tatevsvogel}.

Sections \ref{AB classes} and \ref{Gorenstein dimensions} are a detailed study of the
functors induced on $\bfD(R)$ by those in Theorem~\ref{ifoxby}.  This involves two
different realizations of the derived category as a subcategory of $\bfK(R)$, both
obtained from the localization functor $\bfK(R)\to \bfD(R)$ to $\kproj(R)$: one by
restricting it to the subcategory of K-projective complexes, and the other by restricting
it to $\kinj(R)$, the subcategory of K-injective complexes.  The inclusion $\kproj(R)\to
\bfK(\Proj R)$ admits a right adjoint $\sfp$; for a complex $X$ of projective modules the
morphism $\sfp(X)\to X$ is a K-projective resolution.  In the same way, the inclusion
$\kinj(R)\to \bfK(\Inj R)$ admits a left adjoint $\sfi$, and for a complex $Y$ of
injectives the morphism $Y\to \sfi(Y)$ is a K-injective resolution. 
Consider the functors  $\kauss=\sfi\comp(D\otimes_R-)$ restricted to $\kproj(R)$, and
$\kbass=\sfp\comp\sfq\comp\Hom_R(D,-)$ restricted to $\kinj(R)$.
These functors better visualized as part of the diagram below:
\[
\xymatrixcolsep{2pc} 
\xymatrixrowsep{2.5pc} 
\xymatrix{ 
\bfK(\Proj R) \ar@<1ex>[d]^{\sfp}
  \ar@<-1ex>[rrr]_-{D\otimes_R-}^\sim
  &&& {\bfK(\Inj R)} \ar@<-1ex>[d]_{\sfi} \ar@<-1ex>[lll]_-{\sfq\comp\Hom_R(D,-)} \\
  \kproj(R)\ar@<1ex>[u]^\inc\ar@<-1ex>[rrr]_-{\kauss} &&& {\kinj(R)}
  \ar@<-1ex>[u]_\inc\ar@<-1ex>[lll]_-{\kbass}}
\]
It is clear that $(\kauss,\kbass)$ is an adjoint pair of functors. However, the
equivalence in the upper row of the diagram does not imply an equivalence in the lower
one.  Indeed, given Theorem \ref{ifoxby} and the results in Section \ref{Acyclicity} it is
not hard to prove:

\smallskip

\emph{The natural morphism $X\to \kbass\kauss(X)$ is an isomorphism if and only if the
  mapping cone of the morphism $(D\otimes_RX)\to \sfi(D\otimes_RX)$ is totally acyclic.}

\smallskip

The point of this statement is that the mapping cones of resolutions are, in general, only
acyclic.  Complexes in $\kinj(R)$ for which the morphism $\kauss\kbass(Y)\to Y$ is an
isomorphism can be characterized in a similar fashion; see Propositions~\eqref{aclass} and
\eqref{bclass}.  This is the key observation that allows us to describe, in
Theorems~\eqref{auslander} and \eqref{bass}, the subcategories of $\kproj(R)$ and
$\kinj(R)$ where the functors $G$ and $F$ restrict to equivalences.

Building on these results, and translating to the derived category, we arrive at:

\begin{itheorem}
  A complex $X$ of $R$-modules has finite G-projective dimension if and only if the
  morphism $X\to \RHom_R(D,D\otimes_R^\bfL X)$ in $\bfD(R)$ is an isomorphism and
  $\hh{D\otimes_R^\bfL X}$ is bounded on the left.
\end{itheorem}

The notion of finite G-projective dimension, and finite G-injective dimension, is recalled
in Section \ref{Gorenstein dimensions}. The result above is part of Theorem
\eqref{thm:gproj}; its counterpart for G-injective dimensions is Theorem~\eqref{thm:ginj}.
Given these, it is clear that Theorem I restricts to an equivalence between the category
of complexes of finite G-projective dimension and the category of complexes of finite
G-injective dimension.

Theorems \eqref{thm:gproj} and \eqref{thm:ginj} recover recent results of Christensen,
Frankild, and Holm \cite{Ch:3d}, who arrived at them from a different perspective.  The
approach presented here clarifies the connection between finiteness of G-dimension and
(total) acyclicity, and uncovers a connection between Grothendieck duality and the
equivalence between the categories of complexes of finite G-projective dimension and of
finite G-injective dimension by realizing them as different shadows of the same
equivalence: that given by Theorem \ref{ifoxby}.

So far we have focused on the case where the ring $R$ is commutative.
However, the results carry over, with suitable modifications in the
statements and with nearly identical proofs, to non-commutative rings
that possess dualizing complexes; the appropriate comments are
collected towards the end of each section.  We have chosen to present
the main body of the work, Sections \ref{Foxby
equivalence}--\ref{Gorenstein dimensions}, in the commutative context
in order to keep the underlying ideas transparent, and unobscured by
notational complexity.

\subsection*{Notation}
The following symbols are used to label arrows representing functors
or morphisms: $\sim$ indicates an equivalence (between categories),
$\cong$ an isomorphism (between objects), and $\simeq$ a
quasi-isomorphism (between complexes).

\section{Triangulated categories}
\label{Triangulated categories}

This section is primarily a summary of basic notions and results about
triangulated categories used frequently in this article. For us, the
relevant examples of triangulated categories are homotopy categories
of complexes over noetherian rings; they are the focus of the next
section.  Our basic references are Weibel \cite{Wi}, Neeman
\cite{Ne01}, and Verdier \cite{Ve}.

\begin{bfchunk}{Triangulated categories.}
  Let $\T$ be a triangulated category.  We refer the reader to \cite{Ne01} and \cite{Ve}
  for the axioms that define a triangulated category. When we speak of subcategories, it
  is implicit that they are full.
  
  A non-empty subcategory $\S$ of $\T$ is said to be \emph{thick} if it is a triangulated
  subcategory of $\T$ that is closed under retracts. If, in addition, $\S$ is closed under
  all coproducts allowed in $\T$, then it is \emph{localizing}; if it is closed under all
  products in $\T$ it is \emph{colocalizing}. 

  Let $\C$ be a class of objects in $\T$.  The intersection of the thick subcategories of
  $\T$ containing $\C$ is a thick subcategory, denoted $\Thick(\C)$. We write $\Loc(\C)$,
  respectively, $\Coloc(\C)$, for the intersection of the localizing, respectively,
  colocalizing, subcategories containing $\C$. Note that $\Loc(\C)$ is itself localizing,
  while $\Coloc(\C)$ is colocalizing.
\end{bfchunk}

\begin{bfchunk}{Compact objects and generators.}
\label{compacts}
Let $\T$ be a triangulated category admitting arbitrary coproducts.  An object $X$ of $\T$
is \emph{compact} if $\Hom_{\T}(X,-)$ commutes with coproducts; that is to say, for each
coproduct $\coprod_i Y_i$ of objects in $\T$, the natural morphism of abelian groups
\[
\coprod_{i}\Hom_{\T}(X,Y_i) \lto \Hom_{\T}\bigl(X,\coprod_{i}Y_i\bigr)
\]
is bijective.  The compact objects form a thick subcategory that we denote $\T^c$. We say
that a class of objects $\S$ \emph{generates} $\T$ if $\Loc(\S)=\T$, and that $\T$ is
\emph{compactly generated} if there exists a generating set consisting of compact
objects.

Let $\S$ be a class of compact objects in $\T$. Then $\S$ generates $\T$ if and only if
for any object $Y$ of $\T$, we have $Y=0$ provided that $\Hom_{\T}(\Si^n S,Y)=0$ for all
$S$ in $\S$ and $n\in\bbZ$; see \cite[(2.1)]{Ne96}.
\end{bfchunk}

Adjoint functors play a useful, if technical, role in this work, and pertinent results on
these are collected in the following paragraphs. MacLane's book~\cite[Chapter IV]{Mc} is
the basic reference for this topic; see also \cite[(A.6)]{Wi}.

\begin{bfchunk}{Adjoint functors.}
  Given categories $\A$ and $\B$, a diagram
\[
\xymatrix{ \A \ar@<-1ex>[rr]_-{\sfF}&& \B \ar@<-1ex>[ll]_-{\sfG}}
\]
indicates that $\sfF$ and $\sfG$ are adjoint functors, with $\sfF$ left adjoint to $\sfG$;
that is to say, there is a natural isomorphism $\Hom_{\B}(\sfF(A),B)\cong
\Hom_{\A}(A,\sfG(B))$ for $A\in \A$ and $B\in \B$.
\end{bfchunk}

\begin{chunk}
\label{adjoints:test}
Let $\T$ be a category, $\S$ a full subcategory of $\T$, and $\sfq\col \T\to \S$ a right
adjoint of the inclusion $\inc\colon \S\to\T$.  Then $\sfq\comp\inc \cong {\sf id}_{\S}$.
Moreover, for each $T$ in $\T$, an object $P$ in $\S$ is isomorphic to $\sfq(T)$ if and
only if there is a morphism $P\to T$ with the property that the induced map
$\Hom_{\T}(S,P)\to \Hom_{\T}(S,T)$ is bijective for each $S\in \S$.
\end{chunk}

\begin{chunk}
\label{adjoints}
Let $\sfF\col \S\to \T$ be an exact functor between triangulated categories such that
$\S$ is compactly generated.
\begin{enumerate}[\quad\rm(1)]
\item The functor $\sfF$ admits a right adjoint if and only if it preserves coproducts.
\item The functor $\sfF$ admits a left adjoint if and only if it preserves products.
\item If $\sfF$ admits a right adjoint $\sfG$, then $\sfF$ preserves compactness if and
  only if $\sfG$ preserves coproducts.
\end{enumerate}
For (1), we refer to \cite[(4.1)]{Ne96}; for (2), see \cite[(8.6.1)]{Ne01}; for (3), see
\cite[(5.1)]{Ne96}.
\end{chunk}

\begin{bfchunk}{Orthogonal classes.}
\label{orthogonal classes}
Given a class $\C$ of objects in a triangulated category $\T$, the full subcategories
\begin{alignat*}{2}
  &\phantom{^\perp}\C^\perp & &=\{Y\in\T\mid\Hom_\T(\Si^nX,Y)=0
  \quad \textrm{for all $X\in\C$ and $n\in\bbZ$}\}\,,\\
  &^\perp\C& &=\{X \in\T\mid\Hom_\T(X,\Si^nY)=0 \quad \textrm{for all $Y\in\C$ and
    $n\in\bbZ$}\}\,.
\end{alignat*}
are called the classes \emph{right orthogonal} and \emph{left orthogonal} to $\C$,
respectively.  It is elementary to verify that $\C^\perp$ is a colocalizing subcategory of
$\T$, and equals $\Thick(\C)^\perp$. In the same vein, $^\perp\C$ is a localizing
subcategory of $\T$, and equals ${}^\perp\Thick(\C)$.

Caveat: Our notation for orthogonal classes conflicts with the one in \cite{Ne01}.
\end{bfchunk}

An additive functor $\sfF\col \A\to \B$ between additive categories is
an \emph{equivalence up to direct factors} if $\sfF$ is full and
faithful, and every object in $\B$ is a direct factor of some object
in the image of $\sfF$.
 
\begin{proposition}
\label{pr:quotients}
Let $\T$ be a compactly generated triangulated category and let $\C\subseteq\T$ be a class
of compact objects.
\begin{enumerate}[\quad\rm(1)]
\item The triangulated category $\C^\perp$ is compactly generated. The inclusion
  $\C^\perp\to\T$ admits a left adjoint which induces, up to direct factors, an
  equivalence
\[
\T^c/{\Thick(\C)} \xlto{\sim} (\C^\perp)^c\,.
\]
\item For each class $\B\subseteq\C$, the triangulated category $\B^\perp/\C^\perp$ is
  compactly generated. The canonical functor $\B^\perp\to\B^\perp/\C^\perp$ induces, up to
  direct factors, an equivalence
\[
\Thick(\C)/{\Thick(\B)}\xlto{\sim}(\B^\perp/\C^\perp)^c\,.
\]
\end{enumerate}
\end{proposition}

\begin{proof} 
First observe that $\C$ can be replaced by a set of objects because
  the isomorphism classes of compact objects in $\T$ form a set.
  Neeman gives in \cite[(2.1)]{Ne92} a proof of (1); see also
  \cite[p.~553 ff]{Ne92}. For (2), consider the following diagram
\[
\xymatrix{ \T^c\ar[rr]^-\can\ar[d]^\inc && \T^c/\Thick(\B)\ar[rr]^-\can\ar[d]
  && \T^c/\Thick(\C)\ar[d]\\
  \T\ar@<-1ex>[rr]_-{\sfa} && \B^\perp\ar@<-1ex>[rr]_-{\sfb} \ar@<-1ex>[ll]_-{\inc} &&
  \C^\perp\ar@<-1ex>[ll]_-{\inc} }
\]
where $\sfa$ and $\sfb$ denote adjoints of the corresponding inclusion functors and
unlabeled functors are induced by $\sfa$ and $\sfb$ respectively. The localizing
subcategory $\Loc(\C)$ of $\T$ is generated by $\C$ and hence it is compactly generated
and its full subcategory of compact objects is precisely $\Thick(\C)$; see
\cite[(2.2)]{Ne92}.  Moreover, the composite
\[
\Loc(\C)\xto{\inc}\T\xto{\can}\T/\C^\perp
\]
is an equivalence. From the right hand square one obtains an analogous description of
$\B^\perp/\C^\perp$, namely: the objects of $\C$ in $\T^c/\Thick(\B)$ generate a
localizing subcategory of $\B^\perp$, and this subcategory is compactly generated and
equivalent to $\B^\perp/\C^\perp$. Moreover, the full subcategory of compact objects in
$\B^\perp/\C^\perp$ is equivalent to the thick subcategory generated by $\C$ which is, up
to direct factors, equivalent to $\Thick(\C)/\Thick(\B)$.
\end{proof}

\section{Homotopy categories}
\label{Homotopy categories}

We begin this section with a recapitulation on the homotopy category of an additive
category. Then we introduce the main objects of our study: the homotopy categories of
projective modules, and of injective modules, over a noetherian ring, and establish
results which prepare us for the development in the ensuing sections.

Let $\A$ be an additive category; see \cite[(A.4)]{Wi}.  We grade complexes
cohomologically, thus a complex $X$ over $\A$ is a diagram
\[
\cdots \lto X^n \xlto{\dd^n} X^{n+1}\xlto{\dd^{n+1}} X^{n+2}\lto \cdots
\]
with $X^n$ in $\A$ and $\dd^{n+1}\comp\dd^n=0$ for each integer $n$.  For such a complex
$X$, we write $\Si X$ for its suspension: $(\Si X)^n=X^{n+1}$ and $\dd_{\Si X}= - \dd_X$.

Let $\bfK(\A)$ be the homotopy category of complexes over $\A$; its objects are complexes
over $\A$, and its morphisms are morphisms of complexes modulo homotopy equivalence.  The
category $\bfK(\A)$ has a natural structure of a triangulated category; see \cite{Ve} or
\cite{Wi}.

Let $R$ be a ring. Unless stated otherwise, modules are left modules; right modules are
sometimes referred to as modules over $R^\op$, the opposite ring of $R$. This proclivity
for the left carries over to properties of the ring as well: when we say noetherian
without any further specification, we mean left noetherian, etc.  We write $\bfK(R)$ for
the homotopy category of complexes over $R$; it is $\bfK(\A)$ with $\A$ the category of
$R$-modules.  The paragraphs below contain basic facts on homotopy categories required in
the sequel.

\begin{chunk}
\label{calculus}
Let $\A$ be an additive category, and let $X$ and $Y$ complexes over
$\A$.  Set $\bfK=\bfK(\A)$. Let $d$ be an integer.  We write $X^{\ges
d}$ for the subcomplex $$\cdots\to 0\to X^{d}\to X^{d+1}\to\cdots$$ of $X$, and
$X^{\les d-1}$ for the quotient complex $X/X^{\ges d}$. In $\bfK$
these fit into an exact triangle
\[
X^{\ges d} \lto X \lto X^{\les d-1}\lto \Si X^{\ges d} \tag{$\ast$}
\] 
This induces homomorphisms of abelian groups $\Hom_{\bfK}(X,Y)\to \Hom_{\bfK}(X^{\ges
  d},Y)$ and  $\Hom_{\bfK}(X^{\les d-1},Y)\to \Hom_{\bfK}(X,Y)$.
These have the following properties.
\begin{enumerate}[\quad\rm(1)]
\item One has isomorphisms of abelian groups:
\[
\HH d{\Hom_{\A}(X,Y)}\cong \Hom_{\bfK}(X,\Si^dY)\cong \Hom_{\bfK}(\Si^{-d}X,Y)\,.
\]
\item If $Y^n=0$ for $n\ge d$, then the map \( \Hom_{\bfK}(X^{\les d},Y)\to \Hom_{\bfK}(X,
  Y) \) is bijective.
\item If $Y^n=0$ for $n\le d$, then the map $\Hom_{\bfK}(X,Y)\to \Hom_{\bfK}(X^{\ges d},
  Y)$ is bijective.
\end{enumerate}
There are also versions of (2) and (3), where the hypothesis is on $X$.

Indeed, these remarks are all well-known, but perhaps (2) and (3) less so than (1).  To
verify (2), note that (1) implies
\[
\HH 0{\Hom_{\A}(X^{\ges d+1},Y)}=0=\HH 1{\Hom_{\A}(X^{\ges d+1},Y)}\,,
\]
so applying $\Hom_{\A}(-,Y)$ to the exact triangle ($\ast$) yields that the induced
homomorphism of abelian groups
\[
\HH 0{\Hom_{\A}(X^{\les d},Y)}\lto \HH 0{\Hom_{\A}(X,Y)}
\]
is bijective, which is as desired. The argument for (3) is similar.
\end{chunk}

Now we recall, with proof, a crucial observation from \cite[(2.1)]{Kr}:

\begin{chunk}
\label{calculus:II}
Let $R$ be a ring, $M$ an $R$-module, and let $\sfi M$ be an injective resolution of $M$.
Set $\bfK=\bfK(R)$. If $Y$ is a complex of injective $R$-modules, the induced map
\[
\Hom_{\bfK}(\sfi M,Y)\lto \Hom_{\bfK}(M,Y)
\]
is bijective. In particular, $\Hom_{\bfK}(\sfi R,Y)\cong \HH 0Y$.

Indeed, one may assume $(\sfi M)^n=0$ for $n \leq -1$, since all injective resolutions of
$M$ are isomorphic in $\bfK$.  The inclusion $M\to \sfi M$ leads to an exact sequence of
complexes
\[
0\lto M\lto \sfi M \lto X\lto 0
\]
with $X^n=0$ for $n\leq -1$ and $\hh X=0$.  Therefore for $d=-1,0$ one has isomorphisms
\[
\Hom_{\bfK}(\Si^dX,Y) \cong \Hom_{\bfK}(\Si^dX,Y^{\ges -1}) = 0\,,
\]                      
where the first one holds by an analogue of (\ref{calculus}.2), and the second holds
because $Y^{\ges -1}$ is a complex of injectives bounded on the left.  It now follows from
the exact sequence above that the induced map $\Hom_{\bfK}(\sfi M,Y)\to \Hom_{\bfK}(M,Y)$
is bijective.
\end{chunk}
  
The results below are critical ingredients in many of our arguments. We write
$\bfK^{-,b}(\proj R)$ for the subcategory of $\bfK(R)$ consisting of complexes $X$ of
finitely generated projective modules with $\hh X$ bounded and $X^n=0$ for $n\gg0$, and
$\bfD^f(R)$ for its image in $\bfD(R)$, the derived category of $R$-modules.

\begin{chunk}
\label{compactgen}
Let $R$ be a (not necessarily commutative) ring.

\begin{enumerate}[\quad\rm(1)]
\item When $R$ is coherent on both sides and flat $R$-modules have finite projective
  dimension, the triangulated category $\bfK(\Proj R)$ is compactly generated and the
  functors $\Hom_R(-,R)\col \bfK(\Proj R)\to \bfK(R^\op)$ and $\bfK(R^\op) \to \bfD(R^\op)$ induce
  equivalences
\[
\bfK^c(\Proj R)\xlto\sim \bfK^{-,b}(\proj R^\op)^\op \xlto\sim \bfD^f(R^\op)^\op.
\]
\item When $R$ is  noetherian, the triangulated category $\bfK(\Inj R)$ is compactly
  generated, and the canonical functor $\bfK(\Inj R)\to\bfD(R)$ induces an equivalence
\[
\bfK^c(\Inj R)\xlto\sim \bfD^f(R)
\]
\end{enumerate}  
Indeed, (1) is a result of J{\o}rgensen \cite[(2.4)]{Jo} and (2) is a result of Krause
\cite[(2.3)]{Kr}.
\end{chunk}

In the propositions below $d(R)$ denotes the supremum of the projective dimensions of all
flat $R$-modules.
\begin{proposition}
\label{flat:projection}
Let $R$ be a two-sided coherent ring such that $d(R)$ is finite.  The
inclusion $\bfK(\Proj R)\to \bfK(\Flat R)$ admits a right adjoint:
\[
\xymatrix{ {\bfK(\Proj R)} \ar@<-1ex>[rr]_-{\inc}&& \bfK(\Flat R) \ar@<-1ex>[ll]_-{\sfq}}
\]
Moreover, the category $\bfK(\Proj R)$ admits arbitrary products.
\end{proposition}

\begin{proof}
  By Proposition (\ref{compactgen}.1), the category $\bfK(\Proj R)$ is compactly
  generated.  The inclusion $\inc$ evidently preserves coproducts, so (\ref{adjoints}.1)
  yields the desired right adjoint $\sfq$. The ring $R$ is right coherent, so the
  (set-theoretic) product of flat modules is flat, and furnishes $\bfK(\Flat R)$ with a
  product.  Since $\inc$ is an inclusion, the right adjoint $\sfq$ induces a product on
  $\bfK(\Proj R)$: the product of a set of complexes $\{P_\lambda\}_{\lambda\in\Lambda}$
  in $\bfK(\Proj R)$ is the complex $\sfq\bigl(\prod_{\lambda} P_\lambda\bigr)$.
\end{proof}

The proof of Theorem \ref{pi:form} below uses homotopy limits in the homotopy category of
complexes; its definition is recalled below.

\begin{bfchunk}{Homotopy limits.}
\label{holim}
Let $R$ be a ring and let $\cdots \to X(r+1)\to X(r) $ be a sequence of morphisms in
$\bfK(R)$. The \emph{homotopy limit} of the sequence $\{X(i)\}$, denoted $\holim X(i)$, is
defined by an exact triangle
\[
\xymatrixrowsep{2pc} 
\xymatrixcolsep{1.5pc} \xymatrix{
\holim X(i)\ar@{->}[r]
     &\prod_{i\ges r} X(i) \ar@{->}[rr]^-{\id-\shift}
    &&\prod_{i\ges r} X(i)\ar@{->}[r] & \Si\holim X(i) \,.}
\]
The homotopy limit is uniquely defined, up to an isomorphism in $\bfK(R)$; see \cite{BN}
for details.
\end{bfchunk}

The result below identifies, in some cases, a homotopy limit in the homotopy category with
a limit in the category of complexes.

\begin{lemma}
\label{limits}
Let $R$ be a ring. Consider a sequence of complexes of $R$-modules:
\[
\cdots \lto X(i)\xlto{\eps(i)}X(i-1)\lto\cdots \lto X(r+1)\xlto{\eps(r+1)} X(r)\,.
\]
If for each degree $n$, there exists an integer $s_n$ such that $\eps(i)^n$ is an
isomorphism for $i\geq s_n+1$, then there exists a degree-wise split-exact sequence of
complexes
\[
\xymatrixrowsep{2pc} 
\xymatrixcolsep{1.5pc} \xymatrix{
0\ar@{->}[r] &\varprojlim X(i)\ar@{->}[r] &\prod_i X(i) \ar@{->}[rr]^-{\id-\shift}
    &&\prod_i X(i)\ar@{->}[r]  &0\,.}
\]
In particular, it induces in $\bfK(R)$ an isomorphism $\holim X(i)\cong \varprojlim
X(i)$.
\end{lemma}

\begin{proof}
  To prove the desired degree-wise split exactness of the sequence, it suffices to note
  that if $\cdots \lto M(r+1)\xlto{\delta(r+1)} M(r)$ is a sequence of $R$-modules such
  that $\delta(i)$ is an isomorphism for $i\geq s+1$, for some integer $s$, then one has a
  split exact sequence of $R$-modules:
\[
\xymatrixrowsep{2pc} 
\xymatrixcolsep{1.5pc} \xymatrix{
0\ar@{->}[r] & M(s)\ar@{->}[r]^-{\eta} &\prod_i M(i) \ar@{->}[rr]^-{\id-\shift}
    &&\prod_i M(i)\ar@{->}[r]  &0\,,}
\]
where the morphism $\eta$ is induced by $\eta_i\col M(s)\to M(i)$ with
\[
\eta_i = \begin{cases}
\delta(i+1)\cdots \delta(s) & \text{if $i\leq s-1$}\\
\id &  \text{if $i= s$}\\
\delta(i)^{-1}\cdots \delta(s+1)^{-1} & \text{if $i\geq s+1$}\,.
\end{cases}
\]

Indeed, in the sequence above, the map ($\id-\shift$) is surjective since the
system $\{M_i\}$ evidently satisfies the Mittag-Leffler condition, see \cite[(3.5.7)]{Wi}.
Moreover, a direct calculation shows that $\Im(\eta)=\Ker(\id-\shift)$.  It
remains to note that the morphism $\pi\col \prod M(i)\to M(s)$ defined by $\pi(a_i)=a_s$
is such that $\pi\eta = \id$.

Finally, it is easy to verify that degree-wise split exact sequences of complexes induce
exact triangles in the homotopy category. Thus, by the definition of homotopy limits, see
\eqref{holim}, and the already established part of the lemma, we deduce: $\holim X(i)\cong
\varprojlim X(i)$ in $\bfK(R)$, as desired.
\end{proof}

The result below collects some properties of the functor $\sfq\col\bfK(\Flat
R)\to\bfK(\Proj R)$. It is noteworthy that the proof of part (3) describes an explicit
method for computing the value of $\sfq$ on complexes bounded on the left.  As usual, a
morphism of complexes is called a \emph{\quism} if the induced map in homology is
bijective.

\begin{theorem}
\label{pi:form}
Let $R$ be a two-sided coherent ring with $d(R)$ finite, and let $F$ be a
complex of flat $R$-modules.
\begin{enumerate}[\quad\rm(1)]
\item The morphism $\sfq(F)\to F$ is a \quism.
\item If $F^n=0$ for $n\gg0$, then $\sfq(F)$ is a projective resolution of $F$.
\item If $F^n=0$ for $n\le r$, then $\sfq(F)$ is isomorphic to a complex $P$ with
  $P^{n}=0$ for $n\le r-d(R)$.
\end{enumerate}
\end{theorem}

\begin{proof} (1)
  For each integer $n$, the map $\Hom_\bfK(\Si^nR,\sfq(F))\to\Hom_{\bfK}(\Si^nR,F)$,
  induced by the morphism $\sfq(F)\to F$, is bijective; this is because $R$ is in
  $\bfK(\Proj R)$. Therefore (\ref{calculus}.1) yields $\HH{-n}{\sfq(F)}\cong \HH{-n}F$,
  which proves (1).
  
  (2) When $F^n=0$ for $n\ge r$, one can construct a projective resolution $P\to F$ with
  $P^n=0$ for $n\ge r$. Thus, for each $X\in \bfK(\Proj R)$ one has the diagram below
\[
\Hom_{\bfK}(X^{\les r},P) = \Hom_{\bfK}(X,P) \to \Hom_{\bfK}(X,F)= \Hom_{\bfK}(X^{\les
  r},F)\,.
\]
where equalities hold by (\ref{calculus}.2).  The complex $X^{\les r}$ is K-projective,
so the composed map is an isomorphism; hence the same is true of the one in the middle.
This proves that $\sfq(F)\cong P$; see \eqref{adjoints:test}.

(3) We may assume $d(R)$ is finite. The construction of the complex $P$ takes place in the
category of complexes of $R$-modules. Note that $F^{>i}$ is a subcomplex of $F$ for each
integer $i\ge r$\,; denote $F(i)$ the quotient complex $F/F^{>i}$. One has surjective
morphisms of complexes of $R$-modules
\[
\cdots \lto F(i)\xlto{\eps(i)}F(i-1)\lto\cdots \lto F(r+1) \xlto{\eps(r+1)} F(r) = 0
\]
with $\Ker(\eps(i))=\Si^iF^{i}$. The surjections $F\to F(i)$ are compatible with the
$\eps(i)$, and the induced map $F\to\varprojlim F(i)$ is an isomorphism.  The plan is to
construct a commutative diagram in the category of complexes of $R$-modules
\[
\begin{split}
  \xymatrixrowsep{2pc} \xymatrixcolsep{1.5pc} \xymatrix{ \cdots\ar@{->}[r] &
    P(i)\ar@{->}[r]^-{\delta(i)} \ar@{->}[d]_-{\kappa(i)} & P(i-1)\ar@{->}[r]
    \ar@{->}[d]_-{\kappa(i-1)}
    &\cdots\ar@{->}[r]& P(r+1)\ar@{->}[r]^-{\delta(r+1)}\ar@{->}[d]_-{\kappa(r+1)}& P(r)=0\\
    \cdots\ar@{->}[r] & F(i)\ar@{->}[r]^-{\eps(i)} & F(i-1)\ar@{->}[r] &\cdots
    \ar@{->}[r]& F(r+1) \ar@{->}[r]^-{\eps(r+1)}& F(r)=0}
\end{split}\tag{$\dagger$}
\]
with the following properties: for each integer $i\ge r+1$ one has that
\begin{enumerate}[\quad\rm(a)]
\item $P(i)$ consists of projectives $R$-modules and $P(i)^n=0$ for $n\not\in (r-d(R),i]$;
\item $\delta(i)$ is surjective, and $\Ker\delta(i)^n=0$ for $n< i - d(R)$;
\item $\kappa(i)$ is a surjective \quism.
\end{enumerate}
The complexes $P(i)$ and the attendant morphisms are constructed iteratively, starting
with $\kappa(r+1)\col P(r+1)\to F(r+1)=\Si^{r+1}F^{r+1}$ a surjective projective
resolution, and $\delta(r+1)=0$. One may ensure $P(r+1)^n=0$ for $n\ge r+2$, and also for
$n\le r - d(R)$, because the projective dimension of the flat $R$-module $F^{r+1}$ is at
most $d(R)$.  Note that $P(r+1)$, $\delta(r+1)$, and $\kappa(r+1)$ satisfy conditions
(a)--(c).

Let $i\ge r+2$ be an integer, and let $\kappa(i-1)\col P(i-1)\to F(i-1)$ be a homomorphism
with the desired properties. Build a diagram of solid arrows
\[
\xymatrixrowsep{2pc} \xymatrixcolsep{2pc} \xymatrix{
 0\ar@{->}[r] & Q \ar@{-->}[r]
  \ar@{->}[d]^{\theta} & P(i)\ar@{-->}[r]^-{\delta(i)}\ar@{-->}[d]^-{\kappa(i)}
  & P(i-1)\ar@{->}[r]\ar@{->}[d]^{\kappa(i-1)} & 0  \\
  0\ar@{->}[r] &\Si^{i}F^{i}\ar@{->}[r]^-{\iota} & F(i)\ar@{->}[r]^-{\eps(i)}&
  F(i-1)\ar@{->}[r] & 0 }
\]
where $\iota$ is the canonical injection, and $\theta\col Q \to \Si^{i}F^{i}$ is a
surjective projective resolution, chosen such that $Q^n=0$ for $n< i-d(R)$. The Horseshoe
Lemma now yields a complex $P(i)$, with underlying graded $R$-module $Q\oplus P(i-1)$, and
dotted morphisms that form the commutative diagram above; see \cite[(2.2.8)]{Wi}.  It is
clear that $P(i)$ and $\delta(i)$ satisfy conditions (a) and (b). As to (c): since both
$\theta$ and $\kappa(i-1)$ are surjective \quism s, so is $\kappa(i)$.  This completes the
construction of the diagram ($\dagger$).

Set $P={\varprojlim}{P(i)}$; the limit is taken in the category of complexes.  We claim
that $P$ is a complex of projectives and that $\sfq(F)\cong P$ in $\bfK(\Proj R)$.

Indeed, by property (b), for each integer $n$ the map $P(i+1)^n \to P(i)^n$ is bijective
for $i> n+d(R)$, so $P^n= P(n+d(R))^n$, and hence the $R$-module $P^n$ is projective.
Moreover $P^n=0$ for $n\le r-d(R)$, by (a).  

The sequences of complexes $\{P(i)\}$ and $\{F(i)\}$ satisfy the hypotheses of
Lemma~\eqref{limits}; the former by construction, see property (b), and the latter by
definition.  Thus, Lemma~\eqref{limits} yields the following isomorphisms in $\bfK(R)$:
\[
\holim P(i) \cong P \qquad  \text{and} \qquad \holim F(i) \cong F\,.
\]
Moreover, the $\kappa(i)$ induce a morphism $\kappa\col\holim P(i)\to \holim F(i)$ in
$\bfK(R)$.  Let $X$ be a complex of projective $R$-modules. To complete the proof of (3),
it suffices to prove that for each integer $i$ the induced map
\[
\Hom_{\bfK}(X,\kappa(i)) \col \Hom_{\bfK}(X,P(i)) \lto \Hom_{\bfK}(X,F(i))
\]
is bijective. Then, a standard argument yields that $\Hom_{\bfK}(X,\kappa)$ is bijective,
and in turn this implies $P\cong \holim P(i) \cong \sfq(\holim F(i))\cong\sfq(F)$, see
\eqref{adjoints:test}.

Note that, since $\kappa(i)$ is a \quism\, and $P(i)^n=0=F(i)^n$ for $n\ge i+1$, the
morphism $\kappa(i)\col P(i)\to F(i)$ is a projective resolution. Since projective
resolutions are isomorphic in the homotopy category, it follows from (2) that $P(i)\cong
\sfq(F(i))$, and hence that the map $\Hom_{\bfK}(X,\kappa(i))$ is bijective, as desired.
Thus, (3) is proved.
\end{proof}

\section{Dualizing complexes}
\label{Dualizing complexes}

Let $R$ be a commutative noetherian ring. In this article, a \emph{dualizing complex} for
$R$ is a complex $D$ of $R$-modules with the following properties:
\begin{enumerate}[\quad\rm(a)]
\item the complex $D$ is bounded and consists of injective $R$-modules;
\item the $R$-module $\HH nD$ is finitely generated for each $n$;
\item the canonical map $R\to \Hom_R(D,D)$ is a \quism.
\end{enumerate}
See Hartshorne \cite[Chapter V]{Ha} for basic properties of dualizing complexes.  The
presence of a dualizing complex for $R$ implies that its Krull dimension is finite.  As to
the existence of dualizing complexes: when $R$ is a quotient of a Gorenstein ring $Q$ of
finite Krull dimension, it has a dualizing complex: a suitable representative of the
complex $\RHom_Q(R,Q)$ does the job. On the other hand, Kawasaki \cite{Ka} has proved that
if $R$ has a dualizing complex, then it is a quotient of a Gorenstein ring.

\begin{chunk}
\label{duality:local}
A dualizing complex induces a contravariant equivalence of categories:
\begin{equation*}
\xymatrix{
\bfD^f(R) \ar@<-1ex>[rrr]_-{\Hom_R(-,D)} &&& \bfD^f(R) \ar@<-1ex>[lll]_-{\Hom_R(-,D)} }
\end{equation*}
This property characterizes dualizing complexes: if $C$ is a complex of $R$-modules such
that $\RHom_R(-,C)$ induces a contravariant self-equivalence of $\bfD^f(R)$, then $C$ is
isomorphic in $\bfD(R)$ to a dualizing complex for $R$; see \cite[(V.2)]{Ha}.  Moreover,
if $D$ and $E$ are dualizing complexes for $R$, then $E$ is \quic\, to
$P\otimes_RD$ for some complex $P$ which is locally free of rank one; that is to say, for
each prime ideal $\fp$ in $R$, the complex $P_\fp$ is quasi-isomorphic $\Si^n R_\fp$ for
some integer $n$; see \cite[(V.3)]{Ha}.
\end{chunk}

\begin{remark}
\label{dualizing:pie}  
Let $R$ be a ring with a dualizing complex. Then, as noted above, the Krull dimension of
$R$ is finite, so a result of Gruson and Raynaud \cite[(II.3.2.7)]{GR} yields that the
projective dimension of each flat $R$-module is at most the Krull dimension of $R$.  The
upshot is that Proposition \eqref{flat:projection} yields an adjoint functor
\[
\xymatrix{ {\bfK(\Proj R)} \ar@<-1ex>[rr]_-{\inc}&& \bfK(\Flat R) \ar@<-1ex>[ll]_-{\sfq}}
\]
and this has properties described in Theorem \eqref{pi:form}. In the remainder of the
article, this remark will be used often, and usually without comment.
\end{remark}

In \cite{Ch:3d}, Christensen, Frankild, and Holm have introduced a notion of a dualizing
complex for a pair of, possibly non-commutative, rings:

\begin{bfchunk}{Non-commutative rings.}
\label{ncdualizing complexes}
In what follows $\langle S,R\rangle$ denotes a pair of rings, where
$S$ is left noetherian and $R$ is left coherent and right noetherian.  This
context is more restrictive than that considered in \cite[Section
1]{Ch:3d}, where it is not assumed that $R$ is left coherent.  We make this
additional hypothesis on $R$ in order to invoke (\ref{compactgen}.1).

\begin{subchunk}
\label{ncd:definition}
A \emph{dualizing complex} for the pair $\langle S, R\rangle$ is complex $D$ of $S$-$R$
bimodules with the following properties:
\begin{enumerate}[\quad\rm(a)]
\item $D$ is bounded and each $D^n$ is an $S$-$R$ bimodule that is injective both as an
  $S$-module and as an $R^\op$-module;
\item $\HH nD$ is finitely generated as an $S$-module and as an $R^\op$-module for each $n$;
\item the following canonical maps are \quism s:
\[
R\lto \Hom_S(D,D)\quad \text{and}\quad S\lto \Hom_{R^\op}(D,D)
\]
\end{enumerate}
\end{subchunk}

When $R$ is commutative and $R=S$ this notion of a dualizing complex coincides with the
one recalled in the beginning of this section.  The appendix in \cite{Ch:3d} contains a
detailed comparison with other notions of dualizing complexes in the non-commutative
context.

The result below implies that the conclusion of Remark \eqref{dualizing:pie}: existence of a
functor $\sfq$ with suitable properties, applies also in the situation considered in
\eqref{ncdualizing complexes}.

\begin{proposition}
\label{nc:duality}
Let $D$ be a dualizing complex for the pair of rings $\langle
S,R\rangle$, where $S$ is left noetherian and $R$ is left coherent and right
noetherian.
\begin{enumerate}[\quad\rm(1)]
\item The projective dimension of each flat $R$-module is finite.
\item The complex $D$ induces a contravariant equivalence:
\begin{equation*}
\xymatrix{
\bfD^f(R^\op) \ar@<-1ex>[rrr]_-{\Hom_{R^\op}(-,D)} &&& \bfD^f(S) \ar@<-1ex>[lll]_-{\Hom_S(-,D)} }
\end{equation*}
\end{enumerate}
\end{proposition}
Indeed, (1) is contained in \cite[(1.5)]{Ch:3d}. Moreover, (2) may be proved as in the
commutative case, see \cite[(V.2.1)]{Ha}, so we provide only a

\begin{proof}[Sketch of a proof of {\rm (2)}]
By symmetry, it suffices to prove that for each complex $X$ of right $R$-modules if $\hh
X$ is bounded and finitely generated in each degree, then so is $\hh{\Hom_{R^\op}(X,D)}$, as
an $S$-module, and that the biduality morphism
\[
\theta(X)\col X\lto \Hom_S(\Hom_{R^\op}(X,D),D))
\]
is a \quism.  To begin with, since $\hh X$ is bounded, we may pass to a \quic\, complex
and assume $X$ is itself bounded, in which case the complex $\Hom_{R^\op}(X,D)$, and hence
its homology, is bounded.

For the remainder of the proof, by replacing $X$ by a suitable
projective resolution, we assume that each $X^i$ is a finitely
generated projective module, with $X^i=0$ for $i\gg 0$. In this
case, for any bounded complex $Y$ of $S$-$R$ bimodules, if the
$S$-module $\hh Y$ is finitely generated in each degree, then so is
the $S$-module $\hh{\Hom_{R^\op}(X,Y)}$; this can be proved by an
elementary induction argument, based on the number $$\sup\{i\mid \HH
iY\ne 0\} - \inf\{i\mid \HH iY\ne 0\}\,,$$ keeping in mind that $S$ is
noetherian. Applied with $Y=D$, one obtains that each $\HH
i{\Hom_{R^\op}(X,D)}$ is finitely generated, as desired.

As to the biduality morphism: fix an integer $n$, and pick an integer $d\leq n$ such that
the morphism of complexes
\[
\Hom_S(\Hom_{R^\op}(X^{\ges d},D),D)) \lto \Hom_S(\Hom_{R^\op}(X,D),D))
\]
is bijective in degrees $\geq n-1$; such a $d$ exists because $D$ is bounded. Therefore,
$\HH n{\theta(X)}$ is bijective if and only if $\HH n{\theta(X^{\ges d})}$ is
bijective. Thus, passing to $X^{\ges d}$, we may assume that $X^i=0$ when $|i|\gg 0$. One
has then a commutative diagram of morphisms of complexes
\[
\xymatrix{
X \otimes_R R\ar@{->}[d]_{\cong}\ar@{->}[r]^-{X\otimes_R\theta(R)} &X\otimes_R \Hom_S(D,D) \ar@{->}[d]_{\cong}\\
X\ar@{->}[r]_-\simeq^-{\theta(X)} 
     & \Hom_S(\Hom_{R^\op}(X,D),D)}
\]
The isomorphism on the right holds because $X$ is a finite complex of finitely generated
projectives; for the same reason, since $\theta(R)$ is a \quism, see
(\ref{ncd:definition}.c), so is $X\otimes_R\theta(R)$.  Thus, $\theta(X)$ is a \quism.
This completes the proof.
\end{proof}
\end{bfchunk}

\section{An equivalence of homotopy categories}
\label{Foxby equivalence}

The standing assumption in the rest of this article is that $R$ is a \emph{commutative}
noetherian ring. Towards the end of each section we collect remarks on the extensions of
our results to the non-commutative context described in \eqref{ncdualizing complexes}.

The main theorem in this section is an equivalence between the homotopy categories of
complexes of projectives and complexes of injectives.  As explained in the discussion
following Theorem \ref{ifoxby} in the introduction, it may be viewed as an extension of
the Grothendieck duality theorem, recalled in \eqref{duality:local}.
Theorem~\eqref{foxby} is the basis for most results in this work.

\begin{remark}
  Let $D$ be a dualizing complex for $R$; see Section \ref{Dualizing complexes}.
  
  For any flat module $F$ and injective module $I$, the $R$-module $I\otimes_RF$ is
  injective; this is readily verified using Baer's criterion. Thus, $D\otimes_R-$ is a
  functor between $\bfK(\Proj R)$ and $\bfK(\Inj R)$, and it factors through $\bfK(\Flat
  R)$. If $I$ and $J$ are injective modules, the $R$-module $\Hom_R(I,J)$ is flat, so
  $\Hom_R(D,-)$ defines a functor from $\bfK(\Inj R)$ to $\bfK(\Flat R)$; evidently it is
  right adjoint to $D\otimes_R-\col\bfK(\Flat R)\to\bfK(\Inj R)$.
\end{remark}

Here is the announced equivalence of categories. The existence of
$\sfq$ in the statement below is explained in Remark
\eqref{dualizing:pie}, and the claims implicit in the right hand side
of the diagram are justified by the preceding remark.

\begin{theorem} 
\label{foxby}
Let $R$ be a noetherian ring with a dualizing complex $D$.  The functor $D\otimes_R-\col
\bfK(\Proj R)\to\bfK(\Inj R)$ is an equivalence.  A quasi-inverse is
$\sfq\comp\Hom_R(D,-)$:
\[
\xymatrix{ \bfK(\Proj R)\ar@<-1ex>[rr]_-{\inc}&& \bfK(\Flat R)\ar@<-1ex>[ll]_-{\sfq}
  \ar@<-1ex>[rr]_-{D\otimes_R -}&& \bfK(\Inj R)\ar@<-1ex>[ll]_-{\Hom_R(D,-)}}
\]
where $\sfq$ denotes the right adjoint of the inclusion $\bfK(\Proj R)\to\bfK(\Flat R)$.
\end{theorem}

\begin{chunk}
  The functors that appear in the theorem are everywhere dense in the remainder of this
  article, so it is expedient to abbreviate them: set
\begin{align*}
  &\auss = D\otimes_R-\, \col \bfK(\Proj R)\lto \bfK(\Inj R)\qquad\text{and} \\
  &\bass = \sfq\comp\Hom_R(D,-)\col \bfK(\Inj R)\lto \bfK(\Proj R) \,.
\end{align*}
The notation `$\auss$' should remind one that this functor is given by a tensor product.
The same rule would call for an `$\mathsf H$' to denote the other functor; unfortunately,
this letter is bound to be confounded with an `$H$', so we settle for an `$\bass$'.
\end{chunk}

\begin{proof}
  By construction, $(\inc,\sfq)$ and $(D\otimes_R-,\Hom_R(D,-))$ are adjoint pairs of
  functors.  It follows that their composition $(\auss, \bass)$ is an adjoint pair of
  functors as well.  Thus, it suffices to prove that $\auss$ is an equivalence: this would
  imply that $S$ is its quasi-inverse, and hence also an equivalence.
  
  Both $\bfK(\Proj R)$ and $\bfK(\Inj R)$ are compactly generated, by
  Proposition~\eqref{compactgen}, and $\auss$ preserves coproducts. It follows, using a
  standard argument, that it suffices to verify that $\auss$ induces an equivalence
  $\bfK^c(\Proj R)\to\bfK^c(\Inj R)$. Observe that each complex $P$ of finitely generated
  projective $R$-modules satisfies
\[
\Hom_R(P,D)\cong D\otimes_R\Hom_R(P,R)\,.
\]
Thus one has the following commutative diagram
\[
\xymatrix{ \bfK^{-,b}(\proj R) \ar[d]^\wr\ar[rr]^-{\Hom_R(-,R)}_-\sim
  && \bfK^c(\Proj R)\ar[rr]^{\auss} &&  \bfK^+(\Inj R)\ar[d]^-\wr \\
  \bfD^f(R)\ar[rrrr]_-{\Hom_R(-,D)} &&&&\bfD^+(R) }
\]
By (\ref{compactgen}.2), the equivalence $\bfK^+(\Inj R)\to\bfD^+(R)$ identifies
$\bfK^c(\Inj R)$ with $\bfD^f(R)$, while by \eqref{duality:local}, the functor
$\Hom_R(-,D)$ induces an auto-equivalence of $\bfD^f(R)$. Hence, by the commutative diagram
above, $\auss$ induces an equivalence $\bfK^c(\Proj R)\to\bfK^c(\Inj R)$.  This completes
the proof.
\end{proof}

In the proof above we utilized the fact that $\bfK(\Proj R)$ and $\bfK(\Inj R)$ admit
coproducts compatible with $\auss$. The categories in question also have products; this is
obvious for $\bfK(\Inj R)$, and contained in Proposition~\eqref{flat:projection} for
$\bfK(\Proj R)$. The equivalence of categories established above implies:

\begin{corollary}
\label{products}
The functors $\auss$ and $\bass$ preserve coproducts and products.
\end{corollary}

\begin{remark}
\label{Dstar}
Let $\sfi R$ be an injective resolution of $R$, and set $\dstar=\bass(\sfi R)$.  Injective
resolutions of $R$ are uniquely isomorphic in $\bfK(\Inj R)$, so the complex $\bass(\sfi
R)$ is independent up to isomorphism of the choice of $\sfi R$, so one may speak of
$\dstar$ without referring to $\sfi R$.

\begin{lemma}
  The complex $\dstar$ is isomorphic to the image of $D$ under the composition
\[
\bfD^f(R)\xlto{\sim}\bfK^{-,b}(\proj R)\xto{\Hom_R(-,R)} \bfK(\Proj R)\,.
\]
\end{lemma}

\begin{proof}
  The complex $D$ is bounded and has finitely generated homology modules, so we may choose
  a projective resolution $P$ of $D$ with each $R$-module $P^n$ finitely generated, and
  zero for $n\gg 0$.  In view of Theorem~\eqref{foxby}, it suffices to verify that
  $\auss(\Hom_R(P,R))$ is isomorphic to $\sfi R$. The complex $\auss(\Hom_R(P,R))$, that
  is to say, $D\otimes_R\Hom_R(P,R)$ is isomorphic to the complex $\Hom_R(P,D)$, which
  consists of injective $R$-modules and is bounded on the left. Therefore $\Hom_R(P,D)$ is
  K-injective. Moreover, the composite 
\[
R\lto \Hom_R(D,D)\lto\Hom_R(P,D)
\]
is a quasi-isomorphism, and one obtains that in $\bfK(\Inj R)$ the complex $\Hom_R(P,D)$
is an injective resolution of $R$.
\end{proof}
\end{remark}

The objects in the subcategory $\Thick(\Proj R)$ of $\bfK(\Proj R)$ are exactly the
complexes of finite projective dimension; those in the subcategory $\Thick(\Inj R)$ of
$\bfK(\Inj R)$ are the complexes of finite injective dimension. It is known that the
functor $D\otimes_R-$ induces an equivalence between these categories; see, for instance,
\cite[(1.5)]{AF:fgd}.  The result below may be read as the statement that this equivalence
extends to the full homotopy categories.

\begin{proposition}
\label{pr:AddD}
Let $R$ be a noetherian ring with a dualizing complex $D$.  The equivalence
$\auss\col\bfK(\Proj R)\to\bfK(\Inj R)$ restricts to an equivalence between $\Thick(\Proj
R)$ and $\Thick(\Inj R)$. In particular, $\Thick(\Inj R)$ equals $\Thick(\Add D)$.
\end{proposition}
\begin{proof}
  It suffices to prove that the adjoint pair of functors $(\auss,\bass)$ in
  Theorem~\eqref{foxby} restrict to functors between $\Thick(\Proj R)$ and $\Thick(\Inj
  R)$.
  
  The functor $\auss$ maps $R$ to $D$, which is a bounded complex of injectives and hence
  in $\Thick(\Inj R)$. Therefore $\auss$ maps $\Thick(\Proj R)$ into $\Thick(\Inj R)$.
  
  Conversely, given injective $R$-modules $I$ and $J$, the $R$-module $\Hom_R(I,J)$ is
  flat. Therefore $\Hom_R(D,-)$ maps $\Thick(\Inj R)$ into $\Thick(\Flat R)$, since $D$ is
  a bounded complex of injectives.  By Theorem (\ref{pi:form}.2), for each flat
  $R$-module $F$, the complex $\sfq(F)$ is a projective resolution of $F$.  The projective
  dimension of $F$ is finite since $R$ has a dualizing complex; see \eqref{dualizing:pie}.
  Hence $\sfq$ maps $\Thick(\Flat R)$ to $\Thick(\Proj R)$.
\end{proof}

\begin{bfchunk}{Non-commutative rings.}
\label{ncfoxby}
Consider a pair of rings $\langle S,R\rangle$ as in \eqref{ncdualizing complexes}, with a
dualizing complex $D$. Given Proposition \eqref{nc:duality}, the proof of
Theorem~\eqref{foxby} carries over verbatim to yield:
\begin{Theorem} 
The functor $D\otimes_R-\col \bfK(\Proj R)\to\bfK(\Inj S)$ is an equivalence, and the
functor $\sfq\comp\Hom_S(D,-)$ is a quasi-inverse. \qed
\end{Theorem}

This basic step accomplished, one can readily transcribe the remaining results in this
section, and their proofs, to apply to the pair $\langle S,R\rangle$; it is clear what the
corresponding statements should be.
\end{bfchunk}

% \begin{remark}
% \label{acyclicity:commutative}
% In the commutative case, $\Hom_R(X,D)$ is acyclic if and only if $X$ is acyclic.
% \end{remark}

\section{Acyclicity versus total acyclicity}
\label{Acyclicity}

This section contains various results concerning the classes of (totally) acyclic
complexes of projectives, and of injectives. We start by recalling appropriate
definitions.

\begin{bfchunk}{Acyclic complexes.}
  A complex $X$ of $R$-modules is \emph{acyclic} if $H^nX=0$ for each integer $n$. We
  denote $\bfK_\ac(R)$ the full subcategory of $\bfK(R)$ formed by acyclic complexes of
  $R$-modules. Set
\[
\bfK_\ac(\Proj R)=\bfK(\Proj R)\cap\bfK_\ac(R) \quad\textrm{and} \quad\bfK_\ac(\Inj
R)=\bfK(\Inj R)\cap\bfK_\ac(R)\,.
\]
\end{bfchunk}

Evidently acyclicity is a property intrinsic to the complex under consideration.  Next we
introduce a related notion which depends on a suitable subcategory of $\Mod R$.

\begin{bfchunk}{Total acyclicity.}
\label{tac}
Let $\A$ be an additive category.  A complex $X$ over $\A$ is \emph{totally acyclic} if
for each object $A\in \A$ the following complexes of abelian groups are acyclic.
\[
\Hom_\A(A,X) \quad\text{and}\quad \Hom_\A(X,A)
\]
  We denote by $\bfK_\tac(\A)$ the full subcategory of $\bfK(\A)$ consisting
of totally acyclic complexes.  Specializing to $\A=\Proj R$ and $A=\Inj R$ one gets the
notion of a \emph{totally acyclic complex of projectives} and a \emph{totally acyclic
  complex of injectives}, respectively. 
\end{bfchunk}

Theorems~\eqref{gdefect:proj} and \eqref{gdefect:inj} below describe various properties of
(totally) acyclic complexes.  In what follows, we write $\bfK^c_\ac(\Proj R)$ and
$\bfK^c_\ac(\Inj R)$ for the class of compact objects in $\bfK_\ac(\Proj R)$ and
$\bfK_\ac(\Inj R)$, respectively; in the same way, $\bfK^c_\tac(\Proj R)$ and
$\bfK^c_\tac(\Inj R)$ denote compacts among the corresponding totally acyclic objects.

\begin{theorem}
\label{gdefect:proj}
Let $R$ be a noetherian ring with a dualizing complex $D$.
\begin{enumerate}[\quad\rm(1)]
\item The categories $\bfK_\ac(\Proj R)$ and $\bfK_\tac(\Proj R)$ are compactly generated.
\item The equivalence $\bfD^f(R)\to \bfK^c(\Proj R)^\op $ induces, up to direct factors,
  equivalences
\begin{gather*}
\bfD^f(R)/\Thick(R) \xlto{\sim}   \bfK^c_\ac(\Proj R)^\op  \\
\bfD^f(R)/\Thick(R,D) \xlto{\sim} \bfK^c_\tac(\Proj R)^\op\,.
\end{gather*}
\item The quotient \( \bfK_\ac(\Proj R)/\bfK_\tac(\Proj R) \) is compactly generated, and
  one has, up to direct factors, an equivalence
\[
\Thick(R,D)/\Thick(R) \xlto{\sim}
 \big[\big(\bfK_\ac(\Proj R)/\bfK_\tac(\Proj R)\big)^c\big]^\op\,.
\]
\end{enumerate}
\end{theorem}

The proof of this result, and also of the one below, which is an analogue for complexes of
injectives, is given in \eqref{gdefect:proofs}. It should be noted that, in both cases,
part (1) is not new: for the one above, see the proof of \cite[(1.9)]{Jo2}, and for the
one below, see \cite[(7.3)]{Kr}.

\begin{theorem}
\label{gdefect:inj}
Let $R$ be a noetherian ring with a dualizing complex $D$.
\begin{enumerate}[\quad\rm(1)]
\item The categories $\bfK_\ac(\Inj R)$ and $\bfK_\tac(\Inj R)$ are compactly generated.
\item The equivalence $\bfD^f(R)\to\bfK^c(\Inj R)$ induces, up to direct factors,
  equivalences
\begin{gather*}
\bfD^f(R)/\Thick(R) \xlto{\sim} \bfK^c_\ac(\Inj R)\\
\bfD^f(R)/\Thick(R,D) \xlto{\sim}  \bfK^c_\tac(\Inj R)\,.
\end{gather*}
\item The quotient \( \bfK_\ac(\Inj R)/\bfK_\tac(\Inj R) \) is compactly generated, and we
  have, up to direct factors, an equivalence
\[
\Thick(R,D)/\Thick(R) \xlto{\sim}
\big(\bfK_\ac(\Inj R)/\bfK_\tac(\Inj R)\big)^c\,.
\]
\end{enumerate}
\end{theorem}

Here is one consequence of the preceding results.  In it, one cannot restrict to complexes
(of projectives or of injectives) of finite modules; see the example in Section \ref{An
  example}.

\begin{corollary} 
\label{tac:gorenstein}
Let $R$ be a noetherian ring with a dualizing complex.  The following conditions are
equivalent.
\begin{enumerate}[\quad\rm(a)]
\item The ring $R$ is Gorenstein.
\item Every acyclic complex of projective $R$-modules is totally acyclic.
\item Every acyclic complex of injective $R$-modules is totally acyclic.
\end{enumerate}
\end{corollary}

\begin{proof}
  Theorems (\ref{gdefect:proj}.3) and (\ref{gdefect:inj}.3) imply that (b) and (c) are
  equivalent, and that they hold if and only if $D$ lies in $\Thick(R)$, that is to say,
  if and only if $D$ has finite projective dimension. This last condition is equivalent to
  $R$ being Gorenstein; see \cite[(V.7.1)]{Ha}.
\end{proof}

\begin{remark}
\label{gdefect:cat}
One way to interpret Theorems (\ref{gdefect:proj}.3) and (\ref{gdefect:inj}.3) is that the
category $\Thick(R,D)/\Thick(R)$ measures the failure of the Gorenstein property for $R$.
This invariant of $R$ appears to possess good functorial properties.  For instance, let
$R$ and $S$ be local rings with dualizing complexes $D_R$ and $D_S$, respectively. If a
local homomorphism $ R\to S$ is quasi-Gorenstein, in the sense of Avramov and
Foxby~\cite[Section 7]{AF:fgd}, then tensoring with $S$ induces an equivalence of
categories, up to direct factors:
\[
-\otimes_R^\bfL S \col \Thick(R,D_R)/\Thick(R) \overset{\sim}\lto \Thick(S,D_S)/\Thick(S)
\]
This is a quantitative enhancement of the ascent and descent of the Gorenstein property
along such homomorphisms.
\end{remark}

The notion of total acyclicity has a useful expression in the notation of
\eqref{orthogonal classes}.
\begin{lemma}
\label{tac:equation}
Let $\A$ be an additive category. One has $\bfK_\tac(\A)= \A^\perp\cap {^\perp\A} $, where
$\A$ is identified with complexes concentrated in degree zero.
\end{lemma}

\begin{proof}
  By (\ref{calculus}.1), for each $A$ in $\A$ the complex $\Hom_\A(X,A)$ is acyclic if and
  only if $\Hom_{\bfK(\A)}(X,\Si^nA)=0$ for every integer $n$; in other words, if and only
  if $X$ is in ${}^\perp\A$.  By the same token, $\Hom_\A(A,X)$ is acyclic if and only if
  $X$ is in $\A^\perp$.
\end{proof}

\begin{chunk}
\label{tac:rings}
Let $R$ be a ring. The following identifications hold:
\begin{align*}
  &\bfK_\tac(\Proj R) = \bfK_\ac(\Proj R)\cap {}^\perp(\Proj R)\\
  &\bfK_\tac(\Inj R)= (\Inj R)^\perp \cap \bfK_\ac(\Inj R)\,.
\end{align*}
Indeed, both equalities are due to \eqref{tac:equation}, once it is observed that for any
complex $X$ of $R$-modules, the following conditions are equivalent: $X$ is acyclic;
$\Hom_R(P,X)$ is acyclic for each projective $R$-module $P$; $\Hom_R(X,I)$ is acyclic for
each injective $R$-module $I$.
\end{chunk}

In the presence of a dualizing complex total acyclicity can be tested against a pair of
objects, rather than against the entire class of projectives, or of injectives, as called
for by the definition. This is one of the imports of the result below.  Recall that $\sfi
R$ denotes an injective resolution of $R$, and that $\dstar=\bass(\sfi R)$; see \eqref{Dstar}.

\begin{proposition}\label{pr:tacprojinj}
  Let $R$ be a noetherian ring with a dualizing complex $D$.
\begin{enumerate}[\quad\rm(1)]
\item The functor $\auss$ restricts to an equivalence of $\bfK_\tac(\Proj R)$ with
  $\bfK_\tac(\Inj R)$.
\item $\bfK_\ac(\Proj R) = \{R\}^\perp$ and $\bfK_\tac(\Proj R)= \{R,\dstar\}^\perp$.
\item $\bfK_\ac(\Inj R) = \{\sfi R\}^\perp$ and $\bfK_\tac(\Inj R) = \{\sfi R,D\}^\perp$.
\end{enumerate}
\end{proposition}

\begin{proof} 
  (1) By Proposition~\eqref{pr:AddD}, the equivalence induced by $\auss$ identifies
  $\Thick(\Proj R)$ with $\Thick(\Inj R)$.  This yields the equivalence below:
\begin{multline*}
  \bfK_\tac(\Proj R) = \Thick(\Proj R)^\perp \cap  {}^\perp \Thick(\Proj R) \\
\xto{\sim} \Thick(\Inj R)^\perp \cap{}^\perp\Thick(\Inj R) = \bfK_\tac(\Inj R)
\end{multline*}
The equalities are by Lemma~\eqref{tac:equation}.
  
(3) That $\bfK_\ac(\Inj R)$ equals $\{\sfi R\}^\perp$ follows from (\ref{calculus:II}).
Given this, the claim on $\bfK_\tac(\Inj R)$ is a consequence of \eqref{tac:rings} and the
identifications
\[
\{D\}^\perp=\Thick(\Add D)^\perp=\Thick(\Inj R)^\perp=(\Inj R)^\perp,
\]
where the second one is due to Proposition~\eqref{pr:AddD}.

(2) The equality involving $\bfK_\ac(\Proj R)$ is immediate from (\ref{calculus}.1). Since
$R\otimes_RD\cong D$ and $\dstar\otimes_RD\cong\sfi R$, the second claim follows from (1) and
(3).
\end{proof}

\begin{bfchunk}{Proof of Theorems~\eqref{gdefect:inj} and \eqref{gdefect:proj}.}
\label{gdefect:proofs}
The category $\T=\bfK(\Inj R)$ is compactly generated, the complexes
$\sfi R$ and $D$ are compact, and one has a canonical equivalence
$\T^c\xto{\sim}\bfD^f(R)$; see (\ref{compactgen}.2).  Therefore,
Theorem~\eqref{gdefect:inj} is immediate from Proposition
(\ref{pr:tacprojinj}.3), and Proposition~\eqref{pr:quotients} applied
with $\B=\{\sfi R\}$ and $\C=\{\sfi R,D\}$.

To prove Theorem~\eqref{gdefect:proj}, set $\T=\bfK(\Proj R)$. By
(\ref{compactgen}.1), this category is compactly generated, and in it
$R$ and $\dstar$ are compact; for $\dstar$ one requires also the
identification in \eqref{Dstar}.  Thus, in view of Proposition
(\ref{pr:tacprojinj}.2), Proposition~\eqref{pr:quotients} applied with
$\B=\{R\}$ and $\C=\{R,\dstar\}$ yields that the categories
$\bfK_\ac(\Proj R)$ and $\bfK_\tac(\Proj R)$, and their quotient, are
compactly generated. Furthermore, it provides equivalences
up to direct factors
\begin{gather*}
\bfK^c(\Proj R)/\Thick(R) \xlto{\sim}  \bfK^c_\ac(\Proj R) \\
\bfK^c(\Proj R)/\Thick(R,\dstar) \xlto{\sim} \bfK^c_\tac(\Proj R) \\
\Thick(R,\dstar)/\Thick(R) \xlto{\sim}
  \big(\bfK_\ac(\Proj R)/\bfK_\tac(\Proj R)\big)^c\,.
\end{gather*}
Combining these with the equivalence $ \bfD^f(R) \to \bfK^c(\Proj R)^\op$ in
(\ref{compactgen}.1) yields the desired equivalences. \qed
\end{bfchunk}

\begin{remark}
  Proposition~(\ref{pr:tacprojinj}.3) contains the following result: a complex of
  injectives $X$ is totally acyclic if and only if both $X$ and $\Hom_R(D,X)$ are acyclic.
  We should like to raise the question: if both $\Hom_R(X,D)$ and $\Hom_R(D,X)$ are
  acyclic, is then $X$ acyclic, and hence totally acylic? An equivalent formulation is: if
  $X$ is a complex of projectives and $X$ and $\Hom_R(X,R)$ are acyclic, is then $X$
  totally acyclic?

  In an earlier version of this article, we had claimed an affirmative answer to this
  question, based on a assertion that if $X$ is a complex of $R$-modules such that
  $\Hom_R(X,D)$ is acyclic, then $X$ is acyclic. This assertion is false. Indeed, let $R$
  be a complete local domain, with field of fractions $Q$. A result of
  Jensen~\cite[Theorem 1]{Je} yields $\Ext^i_R(Q,R)=0$ for $i\geq 1$, and it is easy to
  check that $\Hom_R(Q,R)=0$ as well. Thus, $\Hom_R(Q,\sfi R)$ is acyclic. It remains to
  recall that when $R$ is Gorenstein, $\sfi R$ is a dualizing complex for $R$.
\end{remark}

\begin{bfchunk}{Non-commutative rings.}
\label{ncacyclicity}
Theorems \eqref{gdefect:proj} and \eqref{gdefect:inj}, and
Proposition \eqref{pr:tacprojinj}, all carry over, again with suitable modifications in the
statements, to the pair of rings $\langle S,R\rangle$ from \eqref{ncdualizing complexes}.
The analogue of Corollary~\eqref{tac:gorenstein} is especially interesting:

\begin{Corollary} 
The following conditions are equivalent.
\begin{enumerate}[\quad\rm(a)]
\item The projective dimension of $D$ is finite over $R^\op$.
\item The projective dimension of $D$ is finite over $S$.
\item Every acyclic complex of projective $R$-modules is totally acyclic.
\item Every acyclic complex of injective $S$-modules is totally acyclic. \qed
\end{enumerate}
\end{Corollary}
\end{bfchunk}

\section{An example}
\label{An example}
Let $A$ be a commutative noetherian local ring, with maximal ideal $\fm$, and residue
field $k=A/\fm$.  Assume that $\fm^2=0$, and that $\rank_k(\fm)\ge2$.  Observe that $A$ is
\emph{not} Gorenstein; for instance, its socle is $\fm$, and hence of rank at least $2$.
Let $E$ denote the injective hull of the $R$-module $k$; this is a dualizing complex for
$A$.

\begin{proposition}
\label{example}
Set $\bfK=\bfK(\Proj A)$ and let $X$ be a complex of projective $A$-modules.
\begin{enumerate}[\quad\rm(1)]
\item If $X$ is acyclic and the $A$-module $X^d$ is finite for some $d$, then $X\cong 0$
  in $\bfK$.
\item If $X$ is totally acyclic, then $X\cong 0$ in $\bfK$.
\item The cone of the homothety $A\to \Hom_A(P,P)$, where $P$ is a projective resolution
  of $D$, is an acyclic complex of projectives, but it is not totally acyclic.
\item In the derived category of $A$, one has $\Thick(A,D)=\bfD^f(A)$, and hence
\[
\Thick(A,D)/\Thick(A)=\bfD^f(A)/\Thick(A)\,.
\]
\end{enumerate}
\end{proposition}

The proof is given in \eqref{example:proof}. It hinges on some
properties of minimal resolutions over $A$, which we now recall.
Since $A$ is local, each projective $A$-module is free. The Jacobson
radical $\fm$ of $A$ is square-zero, and in particular, nilpotent.
Thus, Nakayama's lemma applies to each $A$-module $M$, hence it has a
projective cover $P\to M$, and hence a minimal projective resolution;
see \cite[Propositions~3 and 15]{E}. Moreover, $\Omega=\Ker(P\to M)$,
the first syzygy of $M$, satisfies $\Omega\subseteq\fm P$, so that
$\fm\Omega\subseteq \fm^2 P =0$, so $\fm\Omega=0$.

\begin{lemma}
\label{example:lemma}
Let $M$ be an $A$-module; set $b=\ell_A(M)$, $c=\ell_A(\Omega)$.
\begin{enumerate}[\quad\rm(1)]
\item If $M$ is finite, then its Poincar\'e series is
\[
P^A_M(t) = b + \frac{ct}{1-et}
\]
In particular, $\beta^A_n(M)$, the $n$th Betti number of $M$, equals $ce^{n-1}$, for
$n\ge1$.
\item If $\Ext^n_A(M,A)=0$ for some $n\ge 2$, then $M$ is free.
\end{enumerate}
\end{lemma}
\begin{proof}
  (1) This is a standard calculation, derived from the exact sequences
\[
0\lto \fm \lto A\lto k\lto 0\quad\text{and}\quad 0\lto\Omega \lto P\lto M\lto 0
\]
The one on the left implies $P^A_k(t)= 1 + etP^A_k(t)$, so $P^A_k(t) = (1-et)^{-1}$, while
the one on the right yields $P^A_M(t)= b + ctP^A_k(t)$, since $\fm\Omega=0$.

(2) If $M$ is not free, then $\Omega\ne0$ and hence has $k$ as a direct summand.  In this
case, since $\Ext^{n-1}_A(\Omega,A)\cong\Ext^n_A(M,A)=0$, one has $\Ext^{n-1}_A(k,A)=0$,
which in turn implies that $A$ is Gorenstein; a contradiction.
\end{proof}

The following test to determine when an acyclic complex is homotopically trivial is surely
known.  Note that it applies to any (commutative) noetherian ring of finite Krull
dimension, and, in particular, to the ring $A$ that is the focus of this section.

\begin{lemma}
\label{ac:triviality}
Let $R$ be a ring whose finitistic global dimension is finite.  An acyclic complex $X$ of
projective $R$-modules is homotopically trivial if and only if for some integer $s$ the
$R$-module $\Coker(X^{s-1}\to X^s)$ is projective.
\end{lemma}

\begin{proof}
  For each integer $n$ set $M(n)=\Coker(X^{n-1}\to X^n)$. It suffices to prove that the
  $R$-module $M(n)$ is projective for each $n$. This is immediate for $n\le s$ because
  $M(s)$ is projective so that the sequence $\cdots\to X^{s-1}\to X^s \to M(s)\to 0$ is
  split exact.
  
  We may now assume that $n \geq s+1$. By hypothesis, there exists an integer $d$ with the
  following property: for any $R$-module $M$, if its projective dimension, $\pd_RM$ is
  finite, then $\pd_RM\leq d$.  It follows from the exact complex
\[
0\lto M(s)\lto X^{s+1}\lto \cdots \lto X^{n+d}\lto M({n+d})\lto 0
\]
that $\pd_RM({n+d})$ is finite. Thus, $\pd_RM({n+d})\leq d$, and another glance at the
exact complex above reveals that $M(n)$ must be projective, as desired.
\end{proof}

Now we are ready for the

\begin{bfchunk}{Proof of Proposition~\eqref{example}.}
\label{example:proof}
In what follows, set $M(s) = \Coker (X^{s-1}\to X^s)$.

(1) Pick an integer $n\ge 1$ with $e^{n-1}\ge \rank_A(X^d) +1$.  Since $X$ is acyclic,
$\Si^{-d-n}\,X^{\leqslant d+n}$ is a free resolution of the $A$-module $M({n+d})$.  Let
$\Omega$ be the first syzygy of $M({n+d})$. One then obtains the first one of the following
equalities:
\[
\rank_A(X^d) \geq \beta^A_{n}(M({n+d})) \geq \ell_A(\Omega) e^{n-1} \ge
\ell_A(\Omega)(\rank_A(X^d) +1)
\]
The second equality is Lemma (\ref{example:lemma}.1) applied to $M({n+d})$ while the last
one is by the choice of $n$. Thus $\ell_A(\Omega)=0$, so $\Omega=0$ and $M({n+d})$ is free.
Now Lemma \eqref{ac:triviality} yields that $X$ is homotopically trivial.

(2) Fix an integer $d$.  Since $\Si^{-d}\,X^{\les d}$ is a projective resolution of
$M(d)$, total acyclicity of $X$ implies that the homology of $\Hom_A(\Si^{-d}X^{\les
  d},A)$ is zero in degrees $\geq 1$, so $\Ext^n_A(M(d),A)=0$ for $n\ge 1$.  Lemma
(\ref{example:lemma}.2) established above implies $M(d)$ is free. Once again,
Lemma~\eqref{ac:triviality} completes the proof.

(3) Suppose that the cone of $A\to \Hom_A(P,P)$ is totally acyclic.  This leads to a
contradiction: (2) implies that the cone is homotopic to zero, so $A\cong \Hom_A(P,P)$ in
$\bfK$.  This entails the first of the following isomorphisms in $\bfK(A)$; the others are
standard.
\begin{align*}
  \Hom_A(k,A) & \cong \Hom_A(k,\Hom_A(P,P)) \\
  & \cong \Hom_A(P\otimes_Ak,P) \\
  & \cong \Hom_k(P\otimes_Ak,\Hom_A(k,P)) \\
  & \cong \Hom_k(P\otimes_Ak,\Hom_A(k,A)\otimes_AP) \\
  & \cong \Hom_k(P\otimes_Ak,\Hom_A(k,A)\otimes_k(k\otimes_AP))
\end{align*}
Passing to homology and computing ranks yields $\hh {k\otimes_AP}\cong k$, and this
implies $D\cong A$. This cannot be for $\rank_k\soc(D)=1$, while $\rank_k\soc(A)=e$ and
$e\ge2$.

(4) Combining Theorem~(\ref{gdefect:proj}.2) and (3) gives the first part.  The second
part then follows from the first. A direct and elementary argument is also available: As
noted above the $A$-module $D$ is not free; thus, the first syzygy module $\Omega$ of $D$
is non-zero, so has $k$ as a direct summand.  Since $\Omega$ is in $\Thick(A,D)$, we
deduce that $k$, and hence every homologically finite complex of $A$-modules, is in
$\Thick(A,D)$.
\end{bfchunk}

\begin{remark}
\label{tatevsvogel}
Let $A$ be the ring introduced at the beginning of this section, and let $X$ and $Y$ be
complexes of $A$-modules.

The Tate cohomology of $X$ and $Y$, in the sense of J{\o}rgensen \cite{Jo2}, is the
homology of the complex $\Hom_A(T,Y)$, where $T$ is a complete projective resolution of
$X$; see \eqref{complete:res}.  By Proposition (\ref{example}.2) any such $T$, being
totally acyclic, is homotopically trivial, so the Tate cohomology modules of $X$ and $Y$
are all zero.  The same is true also of the version of Tate cohomology introduced by
Krause \cite[(7.5)]{Kr} via complete injective resolutions. This is because $A$ has no
non-trivial totally acyclic complexes of injectives either, as can be verified either
directly, or by appeal to Proposition (\ref{pr:tacprojinj}.1).

These contrast drastically with another generalization of Tate cohomology over the ring
$A$, introduced by Vogel and described by Goichot \cite{Go}.  Indeed, Avramov and
Veliche~\cite[(3.3.3)]{AV} prove that for an arbitrary commutative local ring $R$ with residue
field $k$, if the Vogel cohomology with $X=k=Y$ has finite rank even in a \emph{single}
degree, then $R$ is Gorenstein.
\end{remark}

\section{Auslander categories and Bass categories}
\label{AB classes}

Let $R$ be a commutative noetherian ring with a dualizing complex $D$.  We write
$\kproj(R)$ for the subcategory of $\bfK(\Proj R)$ consisting of K-projective complexes,
and $\kinj(R)$ for the subcategory of $\bfK(\Inj R)$ consisting of K-injective
complexes.  This section is motivated by the following considerations: One has adjoint
pairs of functors
\[
{\xymatrix{
\kproj(R)\ar@<-1ex>[r]_-\inc & \bfK(\Proj R) \ar@<-1ex>[l]_-{\sfp}}}
\quad\text{and}\quad 
{\xymatrix{
\bfK(\Inj R) \ar@<-1ex>[r]_-{\sfi}& \kinj(R) \ar@<-1ex>[l]_-{\inc}}}  
\]
and composing these functors with those in Theorem~\eqref{foxby} gives functors
\[
\kauss = (\sfi\comp\auss) \col \kproj(R) \lto \kinj(R) \quad\text{and}\quad 
\kbass = (\sfp\comp \bass) \col\kinj(R)\lto \kproj(R)\,.
\]
These functors fit into the upper half of the picture below:
\[
\xymatrixcolsep{2pc} 
\xymatrixrowsep{2.5pc} 
\xymatrix{ 
\bfK(\Proj R) \ar@<1ex>[d]^{\sfp}  \ar@<-1ex>[rrr]_-{\auss}^-\sim
  &&& {\bfK(\Inj R)} \ar@<-1ex>[d]_{\sfi} \ar@<-1ex>[lll]_-{\bass} \\
  \kproj(R)\ar@<1ex>[u]^-\inc\ar@<-1ex>[rrr]_-{\kauss}
  &&& {\kinj(R)} \ar@<-1ex>[u]_-\inc\ar@<-1ex>[lll]_-{\kbass} \\
  \bfD(R)\ar@{<-}[u]^\wr \ar@<-1ex>[rrr]_-{D\otimes_R^{\bfL}-} &&& {\bfD(R)}
  \ar@{<-}[u]_\wr \ar@<-1ex>[lll]_-{\RHom_R(D,-)}}
\]
The vertical arrows in the lower half are obtained by factoring the canonical functor
$\bfK(\Proj R)\to\bfD(R)$ through $\sfp$, and similarly $\bfK(\Inj R)\to\bfD(R)$ through
$\sfi$.  A straightforward calculation shows that the functors in the last row of the
diagram are induced by those in the middle. Now, while $\auss$ and $\bass$ are
equivalences -- by Theorem~\eqref{foxby} -- the functors $\kauss$ and $\kbass$ need not
be; indeed, they are equivalences if and only if $R$ is Gorenstein; see Corollary
\eqref{gdim:gor} ahead.  The results in this section address the natural:

\begin{Question}
  Identify subcategories of $\kproj(R)$ and $\kinj(R)$ on which
  $\kauss$ and $\kbass$ restrict to equivalences.
\end{Question}

Given the equivalences in the lower square of the diagram an equivalent problem is to
characterize subcategories of $\bfD(R)$ on which the functors $D\otimes^{\bfL}_R-$ and
$\RHom_R(D,-)$ induce equivalences.  This leads us to the following definitions:

\begin{bfchunk}{Auslander category and Bass category.}
\label{AB:class}
Consider the categories
\begin{align*}
  &{\wh\A}(R) = \{X\in \bfD(R)\mid 
       \text{the natural map $X\to \RHom_R(D,D\otimes^{\bfL}_RX)$ is an isomorphism.}\} \\ 
  &{\wh\B}(R) = \{Y\in \bfD(R)\mid 
       \text{the natural map $D\otimes^{\bfL}_R\RHom_R(D,Y)\to Y$ is an isomorphism.}\}
\end{align*}
The notation is intended to be reminiscent of the ones for the \emph{Auslander category}
$\A(R)$ and the \emph{Bass category} $\B(R)$, introduced by Avramov and Foxby
\cite{AF:fgd}, which are the following subcategories of the derived category:
\begin{align*}
  &\A(R) = \{X\in{\wh\A}(R)\mid \text{$X$ and $D\otimes_R^\bfL X$ are homologically
    bounded.}\} \\ &\B(R) = \{Y\in {\wh\B}(R)\mid \text{$Y$ and $\RHom_R(D,Y)$ are
    homologically bounded.}\}
\end{align*}
The definitions are engineered to lead immediately to the following

\begin{proposition}
  The adjoint pair of functors $(\kauss,\kbass)$ restrict to equivalences of categories
  between ${\wh\A}(R)$ and ${\wh\B}(R)$, and between $\A(R)$ and $\B(R)$.\qed
\end{proposition}

In what follows, we identify ${\wh\A}(R)$ and ${\wh\B}(R)$ with the subcategories of
$\kproj(R)$ and $\kinj(R)$ on which $\bass\circ\auss$ and $\auss\circ\bass$, respectively,
restrict to equivalences. The Auslander category and the Bass category are identified with
appropriate subcategories.
\end{bfchunk}

The main task then is describe the complexes in the categories being considered.  In this
section we provide an answer in terms of the categories of K-projectives and K-injectives;
in the next one, it is translated to the derived category.  Propositions~\eqref{aclass}
and \eqref{bclass} below are the first step towards this end.  In them, the \emph{cone} of
a morphism $U\to V$ in a triangulated category refers to an object $W$ obtained by
completing the morphism to an exact triangle: $U\to V\to W\to \Si U$. We may speak of
\emph{the} cone because they exist and are all isomorphic.

\begin{proposition}
\label{aclass}
Let $X$ be a complex of projective $R$-modules. If $X$ is
K-projective, then it is in $\wh\A(R)$ if and only if the cone of the
morphism $\auss(X) \to \sfi\auss(X)$ in $\bfK(\Inj R)$ is totally
acyclic.
\end{proposition}

\begin{Remark}
  The cone in question is always acyclic, because $\auss(X)\to \sfi\auss(X)$ is an
  injective resolution; the issue thus is the difference between acyclicity and total
  acyclicity.
\end{Remark}

\begin{proof}
  Let $\eta\colon\auss(X)\to \sfi\auss (X)$ be a K-injective resolution. In $\bfK(\Proj
  R)$ one has then a commutative diagram
\begin{gather*}
  \xymatrixcolsep{2pc} \xymatrixrowsep{2pc} \xymatrix{
    X \ar@{->}[r]^-\kappa\ar@{->}[d]_\cong &\kbass\kauss(X)\ar@{->}[d]^-\qis\\
    \bass\auss(X)\ar@{->}[r]^-{\bass(\eta)} &\bass\sfi\auss(X)}
\end{gather*}
of adjunction morphisms, where the isomorphism is by Theorem~\eqref{foxby}.  It is clear
from the diagram above that
\begin{align*}
  \text{$X$ is in $\wh\A(R)$} & \iff \text{$\kappa$ is a \quism} \\
  & \iff \text{$\bass(\eta)$ is a \quism}
\end{align*}
It thus remains to prove that the last condition is equivalent to total acyclicity of the
cone of $\eta$. In $\bfK(\Inj R)$ complete $\eta$ to an exact triangle:
\begin{gather*}
  \xymatrix{ \auss(X)\ar@{->}[r]^-{\eta}_-{\qis}&\sfi\auss(X)\ar@{->}[r] & C \ar@{->}[r]
    &\Si\auss(X)}
\end{gather*}
From this triangle one obtains that $\bass(\eta)$ is a \quism\, if and
only if $\bass(C)$ is acyclic.  Now $\bass(C)$ is quasi-isomorphic to
$\Hom_R(D,C)$, see Theorem (\ref{pi:form}.1), and the acyclicity
of $\Hom_R(D,C)$ is equivalent to $C$ being in $\{D\}^\perp$, in
$\bfK(\Inj R)$.  However, $C$ is already acyclic, and hence in $\{\sfi
R\}^\perp$.  Therefore Proposition (\ref{pr:tacprojinj}.3) implies
that $\bass(C)$ is acyclic if and only if $C$ is totally acyclic, as
desired.
\end{proof}

An analogous proof yields:

\begin{proposition}
\label{bclass}
Let $Y$ be a complex of injective $R$-modules. If $Y$ is K-injective,
then it is in $\wh\B(R)$ if and only if the cone of the morphism
$\sfp\bass(Y) \to \bass(Y)$ in $\bfK(\Proj R)$ is totally
acyclic. \qed
\end{proposition}

\begin{corollary}
\label{gdim:gor}
Let $R$ be a noetherian ring with a dualizing complex. The ring $R$ is
Gorenstein if and only if ${\wh\A}(R)=\kproj(R)$, if and only if
${\wh\B}(R)=\kinj(R)$.
\end{corollary}
\begin{proof}
Combine Propositions~\eqref{aclass} and \eqref{bclass} with
Corollary~\eqref{tac:gorenstein}.
\end{proof}

One shortcoming in Propositions~\eqref{aclass} and \eqref{bclass} is
they do not provide a structural description of objects in the
Auslander and Bass categories.  Addressing this issue requires a
notion of complete resolutions.

\begin{bfchunk}{Complete resolutions.}
\label{complete:res}
The subcategory $\bfK_{\tac}(\Proj R)$ of $\bfK(\Proj R)$ is closed under coproducts;
moreover, it is compactly generated, by Theorem~(\ref{gdefect:proj}.1). Thus, the
inclusion $\bfK_{\tac}(\Proj R)\to \bfK(\Proj R)$ admits a right adjoint:
\[
\xymatrix{ \bfK_{\tac}(\Proj R) \ar@<-1ex>[rr]_-{\inc} && \bfK(\Proj R)
  \ar@<-1ex>[ll]_-{\sft}}
\]
For each complex $X$ in $\bfK(\Proj R)$ we call $\sft(X)$ the \emph{complete projective
  resolution} of $X$.  In $\bfK(\Proj R)$, complete the natural morphism $\sft(X)\to X$ to
an exact triangle:
\[
\sft(X)\lto X\lto \sfu(X)\lto \Si\sft(X)
\]
Up to an isomorphism, this triangle depends only on $X$.

Similar considerations show that the inclusion $\bfK_\tac(\Inj R)\to\bfK(\Inj R)$ admits a
left adjoint. We denote it $\sfs$, and for each complex $Y$ of injectives call $\sfs(Y)$
the \emph{complete injective resolution} of $Y$. This leads to an exact triangle in
$\bfK(\Inj R)$:
\[
\sfv(Y)\lto Y\lto \sfs(Y) \lto \Si \sfv(Y)
\]
\end{bfchunk}

Relevant properties of complete resolutions and the corresponding exact triangles are
summed up in the next two result; the arguments are standard, and details are given for
completeness.

\begin{lemma}
\label{complete:proj}
Let $X$ be a complex of projectives $R$-modules.
\begin{enumerate}[\quad\rm(1)]
\item The morphism $X\to \sfu(X)$ is a \quism\, and $\sfu(X)$ is in $\bfK_{\tac}(\Proj R)^\perp$.
\item Any exact triangle $T\to X\to U\to\Si T$ in $\bfK(\Proj R)$ where $T$ is totally
  acyclic and $U$ is in $\bfK_{\tac}(\Proj R)^\perp$ is isomorphic to $\sft(X)\to X\to
  \sfu(X)\to \Si\sft(X)$.
\end{enumerate}
\end{lemma}

\begin{proof}
  (1) By definition, one has an exact triangle
\[
\sft(X)\lto X\lto \sfu(X)\lto \Si \sft(X)\,.
\]
Since the complex $\sft(X)$ is acyclic, the homology long exact sequence arising from this
triangle proves that $X\to \sfu(X)$ is an \quism, as claimed. Moreover, for each totally
acyclic complex $T$ the induced map below is bijective:
\begin{equation*}
\Hom_{\bfK}(T,\sft(X))\lto\Hom_{\bfK}(T,X)\tag{$\dagger$}
\end{equation*}
This holds because $\sft$ is a right adjoint to the inclusion $\bfK_{\tac}(\Proj R)\to
\bfK(\Proj R)$. Since $\sft(-)$ commutes with translations, the morphism $\Si^n \sft(X)\to
\Si ^nX$ coincides with the morphism $\sft(\Si^n X)\to \Si^n X$. Thus, from ($\dagger$) one
deduces that the induced map
\[
\Hom_{\bfK}(T,\sft(\Si^nX))\lto\Hom_{\bfK}(T,\Si^nX)
\]
is bijective for each integer $n$. It is now immediate from the exact triangle above that
$\Hom_{\bfK}(T,\sfu(X))=0$; this settles (1), since $\bfK_{\tac}(\Proj R)$ is stable under
translations.

(2) Given such an exact triangle, the induced map $\Hom_{\bfK}(-,T)\to \Hom_{\bfK}(-,X)$
is bijective on $\bfK_{\tac}(\Proj R)$, since $\Hom_{\bfK}(-,U)$ vanishes on
$\bfK_{\tac}(\Proj R)$.  Thus, there is an isomorphism $\alpha\col T\to \sft(X)$, by
\eqref{adjoints:test}, and one obtains a commutative diagram
\[
\xymatrixcolsep{2pc} \xymatrixrowsep{2pc} \xymatrix{ T \ar@{->}[r] \ar@{->}[d]^{\alpha} &
  X \ar@{->}[r] \ar@{=}[d] & U\ar@{->}[r]\ar@{-->}[d]^{\beta} &\Si T \ar@{->}[r]
  \ar@{->}[d]^{\Si \alpha}&\cdots \\ \sft(X) \ar@{->}[r] & X \ar@{->}[r]&
  \sfu(X)\ar@{->}[r] &\Si\sft(X) \ar@{->}[r] &\cdots }
\]
of morphisms in $\bfK(\Proj R)$. Since the rows are exact triangles, and we are in a
triangulated category, there exists a $\beta$ as above that makes the diagram commute.
Moreover, since $\alpha$ is an isomorphism, so is $\beta$; this is the desired result.
\end{proof}

One has also a version of Lemma~\eqref{complete:proj} for complexes of injectives; proving
it calls for a new ingredient, provided by the next result.  Recall that $\sfi R$ denotes
an injective resolution of $R$ and $\dstar=\bass(\sfi R)$; see \eqref{Dstar}.

\begin{lemma}
\label{orthogonal:lemma}
${}^\perp\bfK_{\tac}(\Inj R) = \Loc(\sfi R,D)$
\end{lemma}

\begin{proof}
  Proposition (\ref{pr:tacprojinj}.3) implies that $\sfi R$ and $D$ are contained in
  ${}^\perp\bfK_{\tac}(\Inj R)$, and hence so is $\Loc(\sfi R,D)$.  To see that the
  reverse inclusion also holds note that $\Loc(\sfi R,D)$ is compactly generated (by $\sfi
  R$ and $D$) and closed under coproducts. Thus, by (\ref{adjoints}.1), the inclusion
  $\Loc(\sfi R,D)\to \bfK(\Inj R)$ admits a right adjoint, say $\sfr$.  Let $X$ be a
  complex of injectives. Complete the canonical morphism $\sfr(X)\to X$ to an exact
  triangle
\[
\sfr(X)\lto X \lto C \lto \Si \sfr(X)
\]
For each integer $n$ the induced map $\Hom_\bfK(-,\Si^n\sfr(X))\to \Hom_{\bfK}(-,\Si^n X)$
is bijective on $\{\sfi R,D\}$, so the exact triangle above yields that $\Hom_{\bfK}(\sfi
R,\Si^n C)=0=\Hom_{\bfK}(D,\Si^nC)$. Therefore, $C$ is totally acyclic, by Proposition
(\ref{pr:tacprojinj}.3). In particular, when $X$ is in ${}^\perp\bfK_{\tac}(\Inj R)$, one
has $\Hom_{\bfK}(X,C)=0$, so the exact triangle above is split, that is to say, $X$ is a
direct summand of $\sfr(X)$, and hence in $\Loc(\sfi R,D)$, as claimed.
\end{proof}

Here is the analogue of Lemma~\eqref{complete:proj} for complexes of injectives; it is a
better result for it provides a structural description of $\sfv(Y)$.

\begin{lemma}
\label{complete:inj}
Let $Y$ be a complex of injective $R$-modules.
\begin{enumerate}[\quad\rm(1)]
\item The morphism $\sfv(Y)\to Y$ is a \quism\, and $\sfv(Y)$ is in $\Loc(\sfi R,D)$.
\item Any exact triangle $V\to X\to T\to\Si V$ in $\bfK(\Inj R)$ where $T$ is
  totally acyclic and $V$ is in $\Loc(\sfi R,D)$ is isomorphic to
  $\sfv(Y)\to Y\to \sfs(Y)\to\Si\sfv(Y)$.
\end{enumerate}
\end{lemma}

\begin{proof}
  An argument akin to the proof of Lemma (\ref{complete:proj}.1) yields that $\sfv(Y)\to
  Y$ is a \quism\, and that $\sfv(Y)$ is in ${}^\perp\bfK_{\tac}(\Inj R)$, which equals
  $\Loc(\sfi R,D)$, by Lemma~\eqref{orthogonal:lemma}.  Given this, the proof of part (2)
  is similar to that of Lemma (\ref{complete:proj}.2).
\end{proof}

Our interest in complete resolutions is due to Theorems~\eqref{bass} and
\eqref{auslander}, which provide one answer to the question raised at the beginning of
this section.

\begin{theorem}
\label{auslander}
Let $R$ be a noetherian ring with a dualizing complex $D$, and let $X$
be a complex of projective $R$-modules. If $X$ is K-projective, then
the following conditions are equivalent.
\begin{enumerate}[\quad\rm(a)]
\item The complex $X$ is in ${\wh\A}(R)$.
\item The complex $\sfu(X)$ is in $\Coloc(\Proj R)$.
\item In $\bfK(\Proj R)$, there exists an exact triangle $T\to X\to
  U\to\Si U$ where $T$ is totally acyclic and $U$ is in $\Coloc(\Proj
  R)$.
\end{enumerate} 
\end{theorem}

\begin{proof}
  Let $\sft(X)\to X\to \sfu(X)\to\Si\sft(X)$ be the exact triangle associated to the
  complete projective resolution of $X$; see \eqref{complete:res}. Let $\eta \col
  \auss(X)\to\sfi\auss(X)$ be a K-injective resolution, and consider the commutative
  diagram
\[
\xymatrixcolsep{2pc} \xymatrixrowsep{2pc} \xymatrix{ \auss(X)\ar@{->}[r]^-{\eta}_-{\qis}
  \ar@{->}[d]_{\qis} & \sfi\auss(X) \ar@{=}[d] \\ 
  \auss\sfu(X)\ar@{->}[r]_-{\qis}^-{\kappa} & \sfi\auss(X) }
\]
arising as follows: the vertical map on the left is a \quism\, because it sits in the
exact triangle with third vertex $\auss\sft(X)$, which is acyclic since $\sft(X)$ is
totally acyclic; see Proposition (\ref{pr:tacprojinj}.1).  Since $\sfi\auss(X)$ is
K-injective, $\eta$ extends to yield $\kappa$, which is a \quism\, because the other
maps in the square are.

Note that the cone of the morphism $\auss(X)\to \auss\sfu(X)$ is $\Si\,\auss\sft(X)$, so
applying the octahedral axiom to the commutative square above gives us an exact triangle
\begin{gather*}
  \xymatrix{ \Si\auss\sft(X) \ar@{->}[r] & \Cone(\eta) \ar@{->}[r] & \Cone(\kappa)
    \ar@{->}[r] &\Si^2\,\auss\sft(X)}
\end{gather*}
where $\Cone(-)$ refers to the cone of the morphism in parenthesis. Since $\sft(X)$ is
totally acyclic, so is $\auss\sft(X)$, by Proposition (\ref{pr:tacprojinj}.1).  Hence the
exact triangle above yields:
\begin{gather*}
  \text{$\Cone(\eta)$ is totally acyclic if and only if $\Cone(\kappa)$ is totally
    acyclic.}  \tag{$\dagger$}
\end{gather*}
This observation is at the heart of the equivalence one has set out to establish.

(a) $\Rightarrow$ (b): Proposition~\eqref{aclass} yields that $\Cone(\eta)$ is totally
acyclic, and hence so is $\Cone(\kappa)$, by ($\dagger$). Consider the exact triangle
\begin{gather*}
  \xymatrix{ \auss\sfu(X)\ar@{->}[r]^-{\kappa}_-{\qis}& \sfi\auss(X) \ar@{->}[r] &
    \Cone(\kappa) \ar@{->}[r] &\Si\auss\sfu(X)}
\end{gather*}
According to Lemma (\ref{complete:proj}.1) the complex $\sfu(X)$ is in $\bfK_\tac(\Proj
R)^\perp$, so Proposition~\eqref{pr:tacprojinj} yields that $\auss\sfu(X)$ is in
$\bfK_{\tac}(\Inj R)^\perp$, and hence the total acyclicity of $\Cone(\kappa)$ implies
\[
\Hom_{\bfK}(\Cone(\kappa),\auss\sfu(X))=0
\]
Thus the triangle above is split exact, and $\auss\sfu(X)$ is a direct summand of
$\sfi\auss(X)$.  Consequently $\auss\sfu(X)$ is in $\Coloc(\Inj R)$, so, by
Theorem~\eqref{foxby} and Corollary~\eqref{products}, one obtains that $\sfu(X)$ is in
$\Coloc(\Proj R)$, as desired.

(b) $\Rightarrow$ (a): By Corollary~\eqref{products}, as $\sfu(X)$ is in $\Coloc(\Proj R)$
the complex $\auss\sfu(X)$ is in $\Coloc(\Inj R)$, that is to say, it is K-injective.
The map $\kappa\col \auss\sfu(X)\to \sfi\auss(X)$, being a \quism\, between
K-injectives, is an isomorphism. Therefore $\Cone(\kappa)\cong 0$ so ($\dagger$) implies
that $\Cone(\eta)$ is totally acyclic. It remains to recall Proposition~\eqref{aclass}.

That (b) implies (c) is patent, and (c) $\Rightarrow$ (b) follows from
Lemma~\eqref{complete:proj}, because $\bfK_{\tac}(\Proj R)^\perp\supseteq \Coloc(\Proj R)$.
The completes the proof of the theorem.
\end{proof}

An analogous argument yields a companion result for complexes of injectives:

 \begin{theorem}
\label{bass}
Let $R$ be a noetherian ring with a dualizing complex $D$, and let $Y$
be a complex of injective $R$-modules. If $Y$ is K-injective, then the
following conditions are equivalent.
\begin{enumerate}[\quad\rm(a)]
\item The complex $Y$ is in ${\wh\B}(R)$.
\item The complex $\sfv(Y)$ is in $\Loc(D)$.
\item In $\bfK(\Inj R)$, there exists an exact triangle $V\to Y\to T\to\Si V$
  where $V$ is in $\Loc(D)$ and $T$ is totally acyclic.\qed
\end{enumerate}
\end{theorem}

Section \ref{Gorenstein dimensions} translates Theorems \eqref{bass} and \eqref{auslander}
to the derived category of $R$.

\begin{bfchunk}{Non-commutative rings.}
\label{ncabclasses}
Consider a pair of rings $\langle S,R\rangle$ with a dualizing complex $D$, defined in
\eqref{ncdualizing complexes}.  As in \eqref{AB:class}, one can define the Auslander
category of $R$ and the Bass category of $S$; these are equivalent via the adjoint pair of
functors $(D\otimes_R-,\sfq\comp\Hom_S(D,-))$. The analogues of Theorems~\eqref{auslander}
and \eqref{bass} extend to the pair $\langle S,R\rangle$, and they describe the complexes
in $\wh\A(R)$ and $\wh\B(S)$.
 \end{bfchunk}

\section{Gorenstein dimensions}
\label{Gorenstein dimensions}

Let $R$ be a commutative noetherian ring, and let $X$ be a complex of $R$-modules.  We say
that $X$ has \emph{finite Gorenstein projective dimension}, or, in short: \emph{finite
  G-projective dimension}, if there exists an exact sequence of complexes of projective
$R$-modules
\[
0\lto U\lto T\lto \sfp X\lto 0
\]
where $T$ is totally acyclic, $\sfp X$ is a K-projective resolution of $X$, and $U^n=0$
for $n\ll0$.
  
Similarly, a complex $Y$ of $R$-modules has \emph{finite G-injective dimension} if there
exists an exact sequence of complexes of injective $R$-modules
\[
0\lto \sfi Y\lto T\lto V\lto 0
\]
where $T$ is totally acyclic, $\sfi Y$ is a K-injective resolution of $Y$, and $V^n=0$
for $n\gg0$.

The preceding definitions are equivalent to the usual ones, in terms of G-projective and
G-injective resolutions; see Veliche \cite{Vl}, and Avramov and Martsinkovsky \cite{AM}.

The theorem below contains a recent result of Christensen, Frankild, and Holm; more
precisely, the equivalence of (i) and (ii) in \cite[(4.1)]{Ch:3d}, albeit in the case when
$R$ is commutative; however, see \eqref{ncgdim}.

\begin{theorem}
\label{thm:gproj}
Let $R$ be a noetherian ring with a dualizing complex $D$, and $X$ a complex of
$R$-modules.  The following conditions are equivalent:
\begin{enumerate}[\quad\rm(a)]
\item $X$ has finite G-projective dimension.
\item $\sfp X$ is in ${\wh\A}(R)$ and $D\otimes_R^\bfL X$ is homologically
bounded on the left.
\item $\sfu(\sfp X)$ is isomorphic, in $\bfK(\Proj R)$, to a complex $U$
  with $U^n=0$ for $n\ll0$.
\end{enumerate}
When $\hh X$ is bounded, these conditions are equivalent to: $X$ is in $\A(R)$.
\end{theorem}

\begin{proof}
  Substituting $X$ with $\sfp X$, one may assume that $X$ is K-projective and that
  $D\otimes_R^\bfL X$ is \quic\, to $D\otimes_RX$, that is to say, to $\auss(X)$.
  
  \smallskip
  
  (a) $\Rightarrow$ (b): By definition, there is an exact sequence of complexes of
  projectives $0\to U\to T\to X\to 0$ where $T$ is totally acyclic and $U^n=0$ for $n\ll
  0$. Passing to $\bfK(\Proj R)$ gives rise to an exact triangle
\[
U \lto T\lto X \lto \Si U
\]
Since $T$ is totally acyclic, $\auss(X)$ is \quic\, to $\auss(\Si U)$; the latter is
bounded on the left as a complex, hence the former is homologically bounded on the left,
as claimed. This last conclusion yields also that $\auss(\Si U)$ is in $\Coloc(\Inj R)$.
Thus, by Theorem~\eqref{foxby} and Corollary~\eqref{products}, the complex $\Si U$ is in
$\Coloc(\Proj R)$, so the exact triangle above and Theorem \eqref{auslander} imply that
$X$ is in ${\wh\A}(R)$.

\smallskip

(b) $\Rightarrow$ (c): By Theorem~\eqref{auslander}, there is an exact triangle
\[
T \lto X \lto U \lto \Si T
\]
with $T$ totally acyclic and $U$ in $\Coloc(\Proj R)$. The first
condition implies that $\auss(U)$ is \quic\, to $\auss(X)$, and hence
homologically bounded on the left, while the second implies, thanks to
Corollary~\eqref{products}, that it is in $\Coloc(\Inj R)$, that is to
say, it is K-injective.  Consequently $\auss(U)$ is isomorphic to a
complex of injectives $I$ with $I^n=0$ for $n\ll 0$. This implies that
the complex of flat $R$-modules $\Hom_R(D,\auss(U))$ is bounded on the
left.  Theorem (\ref{pi:form}.3) now yields that the complex
$\sfq(\Hom_R(D,\auss(U)))$, that is to say, $\bass\auss(U)$, is
bounded on the left; thus, the same is true of $U$ as it is isomorphic
to $\bass\auss(U)$, by Theorem~\eqref{foxby}. It remains to note that
$\Coloc(\Proj R)\subseteq \bfK_\tac(\Proj R)^\perp$, so $\sfu(X)\cong
U$ by Lemma \eqref{complete:proj}.

\smallskip

(c) $\Rightarrow$ (a): Lift the morphism $X\to \sfu(X)\cong U$ in $\bfK(\Proj R)$ to a
morphism $\alpha\col X\to U$ of complexes of $R$-modules. In the mapping cone exact
sequence
\[
0\lto U\lto \Cone(\alpha)\lto \Si X\lto 0
\]
$\Cone(\alpha)$ is homotopic to $\sft(X)$, and hence totally acyclic, while $U^n=0$ for
$n\ll 0$, by hypothesis. Thus, the G-projective dimension of $\Si X$, and hence of $X$, is
finite.

Finally, when $\hh X$ is bounded, $D\otimes_R^\bfL X$ is always bounded on the right. It
is now clear from definitions that the condition that $X$ is in $\A(R)$ is equivalent to
(b).
\end{proof}

Here is a characterization of complexes in $\bfD(R)$ that are in the Bass category.  For
commutative rings, it recovers \cite[(4.4)]{Ch:3d}; see \eqref{ncgdim}.  The
basic idea of the proof is akin the one for the theorem above, but the details are
dissimilar enough to warrant exposition.

\begin{theorem}
\label{thm:ginj}
Let $R$ be a noetherian ring with a dualizing complex $D$, and $Y$ a
complex of $R$-modules.  The following conditions are equivalent:
\begin{enumerate}[\quad\rm(a)]
\item $Y$ has finite G-injective dimension.
\item $\sfi Y$ is in ${\wh\B}(R)$ and $\RHom_R(D,Y)$ is homologically
bounded on the right.
\item $\sfv(\sfi Y)$ is isomorphic, in $\bfK(\Inj R)$, to a complex $V$
  with $V^n=0$ for $n\gg0$.
\end{enumerate}
When $\hh Y$ is bounded, these conditions are equivalent to: $Y$ is in $\B(R)$.
\end{theorem}

\begin{proof}
  Replacing $Y$ with $\sfi Y$ we assume that $Y$ is K-injective, so $\RHom_R(D,Y)$ is
  \quic\, to $\Hom_R(D,Y)$. In the argument below the following remark is used without
  comment: in $\bfK(\Inj R)$, given an exact triangle
\[
Y_1\lto Y_2\lto T\lto\Si Y_1
\]
if $T$ is totally acyclic, then one has a sequence
\[
\Hom_R(D,Y_1)\xleftarrow{\qis} \bass(Y_1) \xto{\qis} \bass(Y_2)\xto{\qis} \Hom_R(D,Y_2)\,.
\]
of quasi-isomorphisms.  Indeed, the first and the last \quism\, hold by Theorem
(\ref{pi:form}.1), while the middle one holds because $\bass(T)$ is totally acyclic, by
Theorem~\eqref{foxby}.

\smallskip
  
(a) $\Rightarrow$ (b): The defining property of complexes of finite G-injective dimension
provides an exact sequence of complexes of injectives $0\to Y\to T\to V\to 0$ where $T$ is
totally acyclic and $V^n=0$ for $n\gg 0$. Passing to $\bfK(\Inj R)$ gives rise to an exact
triangle
\[
\Si^{-1} V\lto Y\lto T \lto V
\]
Since $T$ is totally acyclic, $\Hom_R(D,\Si^{-1}V)$ is \quic\, to $\Hom_R(D,Y)$; the
former is bounded on the right as a complex, so the latter is homologically bounded on
the right, as claimed. Furthermore, since $V$ is bounded on the right, so is
$\Hom_R(D,\Si^{-1}V)$. Theorem (\ref{pi:form}.2) then yields that $\bass(\Si^{-1}V)$
is its projective resolution, and hence it is in $\Loc(R)$. Thus, by
Theorem~\eqref{foxby}, the complex $\Si^{-1}V$ is in $\Loc(D)$, so the exact triangle
above and Theorem \eqref{bass} imply that $Y$ is in ${\wh\B}(R)$.

\smallskip

(b) $\Rightarrow$ (c): By hypothesis and Theorem~\eqref{bass} there exists and exact
triangle
\[
V \lto Y \lto T\lto\Si V
\]
in $\bfK(\Inj R)$, where $V$ lies in $\Loc(D)$ and $T$ is totally acyclic. Thus $\bass(V)$
is in $\Loc(R)$, that is to say, it is K-projective, and it is \quic\, to $\Hom_R(D,Y)$,
and hence it is homologically bounded on the right. Therefore, $\bass(V)$ is isomorphic to
a complex of projectives $P$ with $P^n=0$ for $n\gg 0$. By Theorem~\eqref{foxby}, this
implies that $V$ is isomorphic to $\auss(P)$, which is bounded on the right.

\smallskip

(c) $\Rightarrow$ (a): Lift the morphism $V\cong \sfv(Y)\to Y$ in $\bfK(\Inj R)$ to a
morphism $\alpha\col V\to Y$ of complexes of $R$-modules. In the mapping cone exact
sequence
\[
0\lto Y\lto \Cone(\alpha)\lto \Si V\lto 0
\]
the complex $\Cone(\alpha)$ is homotopic to $\sfs(Y)$, and hence totally acyclic, while
$V^n=0$ for $n\gg 0$, by hypothesis. Thus, the G-injective dimension of $Y$ is finite.

Finally, when $Y$ is homologically bounded, $\RHom_R(D,Y)$ is bounded on the left, so $Y$
is in $\B(R)$ if and only if it satisfies condition (b).
\end{proof}

\begin{bfchunk}{Non-commutative rings.}
\label{ncgdim}  
Following the thread in \eqref{ncdualizing complexes}, \eqref{ncfoxby},
\eqref{ncacyclicity}, and \eqref{ncabclasses}, the development of this section also
carries over to the context of a pair of rings $\langle S,R\rangle$ with a dualizing
complex $D$. In this case, the analogues of Theorems \eqref{thm:gproj} and
\eqref{thm:ginj} identify complexes of finite G-projective dimension over $R$ and of
finite G-injective dimension over $S$ as those in the Auslander category of $R$ and the
Bass category of $S$, respectively. These results contain \cite[(4.1),(4.4)]{Ch:3d}, but
only when one assumes that the ring $R$ is left coherent as well; the reason for this has
already been given in \eqref{ncdualizing complexes}.
\end{bfchunk}

\section*{Acknowledgments}
This project was initiated in June 2004, when Srikanth Iyengar was visiting the University
of Paderborn. He thanks the Department of Mathematics at Paderborn and its algebra group
for hospitality. The authors thank Lucho Avramov for his comments and suggestions, Lars
W. Christensen for numerous discussions on \cite{Ch:3d} which motivates the results in
Section \ref{Gorenstein dimensions}, and Amnon Neeman for a careful reading of the
manuscript.

\end{document}